\documentclass[a4paper,11pt]{article} 
\usepackage{amsmath,amsfonts} 
\newtheorem{thm}{Theorem}[section]   
\newtheorem{cor}[thm]{Corollary}   
\newtheorem{prop}[thm]{Proposition}   
\newtheorem{lm}[thm]{Lemma}

\newtheorem{rem}[thm]{Remark}  
\newcommand{\RR}{\mathbb{R}}   
\newcommand{\CC}{\mathbb{C}}   
\newcommand{\DD}{\mathbb{D}}   
\newcommand{\TT}{\mathbb{T}}   
\newcommand{\ZZ}{\mathbb{Z}}    
\newcommand{\NN}{\mathbb{N}}

\newcommand{\PP}{\mathbb{P}}
   
\renewcommand{\div}{\mathop{\mathrm{div}}}
\newcommand{\curl}{\mathop{\mathrm{curl}}}
\newcommand{\supp}{\mathop{\mathrm{supp}}}

\newcommand{\dd}{\,\mathrm{d}}
\newcommand{\dD}{\mathrm{d}}
\newcommand{\loc}{\mathrm{loc}}
\newcommand{\cO}{\mathcal{O}}
\newcommand{\one}{\mathbf{1}}
\renewcommand{\epsilon}{\varepsilon}

\begin{document}   
   
\title{\textbf{Global Existence and Long-Time Asymptotics}\\ 
\textbf{for Rotating Fluids in a 3D Layer}}
   
\author{
{\bf Thierry Gallay} \\ 
Institut Fourier (UMR CNRS 5582) \\ 
Universit\'e de Grenoble I \\
B.P. 74 \\
38402 Saint-Martin-d'H\`eres, France \\ 
\texttt{thierry.gallay@ujf-grenoble.fr}\\ 
\\[2mm]
{\bf Violaine Roussier-Michon} \\ 
Institut de Math\'ematiques de Toulouse (UMR CNRS 5219)\\ 
INSA Toulouse \\
135 av Rangueil \\
31077 Toulouse cedex 4, France\\ 
\texttt{roussier@math.univ-toulouse.fr}}
   
\date{December 5, 2008}
   
\maketitle   
\begin{abstract}
The Navier-Stokes-Coriolis system is a simple model for
rotating fluids, which allows to study the influence of 
the Coriolis force on the dynamics of three-dimensional 
flows. In this paper, we consider the NSC system in an infinite
three-dimensional layer delimited by two horizontal planes, 
with periodic boundary conditions in the vertical direction.
If the angular velocity parameter is sufficiently large, 
depending on the initial data, we prove the existence of 
global, infinite-energy solutions with nonzero circulation
number. We also show that these solutions converge toward  
two-dimensional Lamb-Oseen vortices as $t \to \infty$. 
\end{abstract}

\vfill\eject

\section{Introduction}    

In recent years a lot of activity has been devoted to the mathematical
study of geophysical flows, and in particular to various models of
rotating fluids. Taking advantage of the stratification effect due to
the Coriolis force, significant results have been obtained which are
still out of reach for the usual Navier-Stokes system, such as global
existence of solutions for large initial data \cite{bmn,cdgg2} and
stability of boundary layers for small viscosities \cite{grenier,
masmoudi}. We refer the interested reader to the recent monograph 
\cite{cdgg} which contains a general introduction to geophysical
flows, an overview of the mathematical theory, and an extensive 
bibliography.

In this article we study the so-called Navier-Stokes-Coriolis (NSC)
system in a three-dimen\-sional layer delimited by two infinite
horizontal planes, assuming as usual that the rotation vector is
constant and aligned with the vertical axis. This is a reasonably
simple model for the motion of the ocean in a small geographic zone at
mid-latitude, where the variation of the Coriolis force due to the
curvature of Earth can be neglected. More realistic systems exist which
take into account the variations of temperature and salinity inside
the ocean, and include boundary effects modelling the influence of
coasts, the topography of the bottom, or the action of the wind at the
free surface, see \cite{greenspan,pedlosky}. Nevertheless, keeping
only the Coriolis force is meaningful in a first approximation,
because its effect is very important on the ocean's motion at a global
scale due to the fast rotation of Earth compared to typical 
velocities in the ocean.

Our main goal is to investigate the long-time behavior of the
solutions to the NSC system for a fixed, but typically large, value of
the rotation speed. As in \cite{bmn,cdgg2} we shall use the effect of
the Coriolis force to prove global existence of solutions for
large initial data, but the long-time asymptotics of those solutions
turn out to be essentially two-dimensional and are therefore not
affected by the rotation. Thus we shall recover as a leading term in
our expansion the Lamb-Oseen vortex which plays a similar role for the
usual Navier-Stokes system in the plane $\RR^2$ \cite{thgcew2} or the
three-dimensional layer $\RR^2 \times (0,1)$ \cite{roussier}. To avoid
all problems related to boundary layers, we shall always assume that
the fluid motion is {\em periodic} in the vertical direction. This
hypothesis has no physical justification and is only a convenient 
mathematical way to disregard the influence of the boundaries.
Although boundary conditions do play an important role in 
the problem we study and will have to be considered ultimately, 
in this paper we chose to focus on the motion of the fluid in the
bulk. 

We thus consider the Navier-Stokes-Coriolis system in the
three-dimensional layer $\DD = \RR^2 \times \TT^1$, where $\TT^1 =
\RR/\ZZ$ is the one-dimensional torus. The points of $\DD$ will be
denoted by $(x,z)$, where $x = (x_1,x_2) \in \RR^2$ is the horizontal
variable and $z \in \TT^1$ is the vertical coordinate.  The system
reads
\begin{equation}\label{NSC}
  \partial_t u + (u\cdot\nabla)u + \Omega e_3 \wedge u 
  \,=\, \Delta u -\nabla p~, \qquad \div u \,=\,0~, 
\end{equation}
where $u = u(t,x,z) \in \RR^3$ is the velocity field of the fluid, and
$p = p(t,x,z) \in \RR$ is the pressure field. Here and in what
follows, it is understood that differential operators such as $\nabla$
or $\Delta$ act on all spatial variables $(x,z)$, unless otherwise
indicated. System~\eqref{NSC} differs from the usual incompressible
Navier-Stokes equations by the presence of the Coriolis term $\Omega
e_3 \wedge u$, where $\Omega \in \RR$ is a parameter and $e_3 =
(0,0,1)^t$ is the unit vector in the vertical direction. This term is
due to the fact that our reference frame rotates with constant angular
velocity $\Omega/2$ around the vertical axis. Note that \eqref{NSC}
does not contain any centrifugal force, because this effect can be
included in the pressure term $-\nabla p$. For simplicity, the
kinematic viscosity of the fluid has been rescaled to $1$, and the
fluid density has been incorporated in the definition of the pressure
$p$.

As in the ordinary Navier-Stokes system, the role of the pressure
in \eqref{NSC} is to enforce the incompressibility condition 
$\div u = 0$. To eliminate the pressure, one can apply to both 
sides the Leray projector $\PP$, which is just the orthogonal 
projector in $L^2(\DD)^3$ onto the space of divergence-free 
vector fields. This operator has a rather simple expression 
in Fourier variables, which will be given in Appendix~\ref{BSlaw}. 
The projected equation then reads:
\begin{equation}\label{nsft}
  \partial_t u + \PP((u \cdot \nabla)u) + \Omega \PP(e_3\wedge u) 
  \,=\, \Delta u~, \qquad \div u \,=\, 0~.
\end{equation}
Another possibility is to consider the {\em vorticity field} 
$\omega = \curl u$, which satisfies the following evolution 
equation:
\begin{equation}\label{nsvor}
  \partial_t \omega + (u \cdot \nabla)\omega - (\omega \cdot \nabla)u
  - \Omega \partial_z u \,=\, \Delta \omega~.
\end{equation}
Due to the incompressibility condition, the velocity field $u$ can be
reconstructed from the vorticity $\omega$ using the Biot-Savart law,
which in the domain $\DD$ has also a simple expression, see
Appendix~\ref{BSlaw}.

As is clear from \eqref{nsvor}, the vertical coordinate $z$ 
plays a distinguished role in our problem because the rotation 
acts trivially on $z$-independent velocity fields. As a matter 
of fact, even if rotation is absent, the linear evolution 
$\partial_t u = \Delta u$ leads to an exponential decay of 
the fluctuations of $u$ in the vertical direction, due to the
Poincar\'e inequality. For these reasons, it is appropriate to 
decompose the velocity field as $u(t,x,z) = \bar u(t,x) + \tilde 
u(t,x,z)$, where
\begin{equation}\label{Qdef}
  \bar u(t,x) \,=\, (Qu)(t,x) \,\equiv\, \int_{\TT^1} u(t,x,z)
  \dd z
\end{equation}
is the average of $u$ with respect to the vertical variable, and the
remainder $\tilde u = (1-Q)u$ has zero vertical average. We shall say
that $\bar u$ is a {\em two-dimensional} vector field in the sense
that it depends only on the spatial variable $x \in \RR^2$, not on
$z$, but one should keep in mind that $\bar u$ is not necessarily {\em
  horizontal} because its third component $\bar u_3$ is usually
nonzero. A similar decomposition holds for the vorticity, and it is
easy to verify that $\bar \omega = \curl \bar u$ and $\tilde \omega =
\curl \tilde u$. In particular, since $\partial_1 \bar u_1 +
\partial_2 \bar u_2 = 0$ and $\partial_1 \bar u_2 - \partial_2 \bar
u_1 = \bar \omega_3$, the horizontal part of the two-dimensional
velocity field $\bar u$ can be reconstructed from the third component
of the vorticity $\bar \omega$ via the two-dimensional Biot-Savart
law, see Appendix~\ref{BSlaw}. This means that the averaged velocity
field $\bar u(t,x)$ can be represented by two scalar quantities,
namely $\bar u_3(t,x)$ and $\bar \omega_3(t,x)$.

We shall solve the Cauchy problem for equation \eqref{nsft} in 
the Banach space $X$ defined by
\begin{equation}\label{Xdef}
  X \,=\, \Bigl\{u \in H^1_\loc(\DD)^3\,\Big|\, \div u = 0\,,~
  \tilde u \in H^1(\DD)^3\,,~ \bar u_3 \in H^1(\RR^2)\,,~
  \bar \omega_3 \in L^1(\RR^2) \cap L^2(\RR^2)\Bigr\}~,
\end{equation}
equipped with the norm
\[
  \|u\|_X \,=\, \|\tilde u\|_{H^1(\DD)} + \|\bar u_3\|_{H^1(\RR^2)}
  + \|\bar \omega_3\|_{L^1(\RR^2)} + \|\bar \omega_3\|_{L^2(\RR^2)}~.
\]
Observe that $X \not\subset H^1(\DD)^3$, because the two-dimensional
horizontal velocity field $\bar u_h = (\bar u_1,\bar u_2)$ is not
assumed to be square integrable. This slightly unusual choice is
motivated by our desire to include {\em infinite-energy} solutions,
which play a crucial role in the long-time asymptotics of the
Navier-Stokes equations \cite{thgcew,thgcew2}. The most important
example of such a solution is the {\em Lamb-Oseen} vortex, whose
velocity and vorticity fields are given by the following expressions:
\begin{align}\label{Oseen1}
  u^G(t,x) \,&=\, \frac{1}{\sqrt{1+t}}\,U^G\Bigl(\frac{x}{
  \sqrt{1+t}}\Bigr)~, \quad \hbox{where}\quad 
  U^G(\xi) \,=\, \frac{1-e^{-|\xi|^2/4}}{2\pi|\xi|^2}\,
  \begin{pmatrix} -\xi_2 \\ \xi_1 \\ 0 \end{pmatrix}~,\\ \label{Oseen2}
  \omega^G(t,x) \,&=\, \frac{1}{1+t}\,G\Bigl(\frac{x}{
  \sqrt{1+t}}\Bigr)~, \qquad~ \hbox{where}\quad 
  G(\xi) \,=\, \frac{1}{4\pi}\,e^{-|\xi|^2/4}
  \begin{pmatrix} \,0\, \\ 0 \\ 1 \end{pmatrix}~.
\end{align}
As is easily verified, for any $\alpha \in \RR$ and any $\Omega \in
\RR$, the vortex $u(t,x,z) \,=\, \alpha u^G(t,x)$ is an exact 
solution of the NSC system \eqref{nsft}. In fact, one 
has $\PP(u^G\cdot\nabla)u^G = 0$ and $\PP(e_3\wedge u^G) = 0$, so that
$u^G$ solves the linear heat equation $\partial_t u = \Delta u$. 

\medskip
We are now in position to formulate our main result:

\begin{thm}
\label{main}
For any initial data $u_0 \in X$, there exists $\Omega_0 \ge 0$ 
such that, for all $\Omega \in \RR$ with $|\Omega| \ge \Omega_0$, 
the NSC system \eqref{nsft} has a unique global (mild) solution
$u \in C^0([0,\infty),X)$ satisfying $u(0) = u_0$. Moreover
$\|u(t,\cdot) - \alpha u^G(t,\cdot)\|_X \to 0$ as $t \to \infty$, 
where
\begin{equation}\label{alphadef}
  \alpha \,=\, \int_{\DD} (\curl u_0)_3 \dd x \dd z~.
\end{equation}
\end{thm}

This theorem contains in fact two different statements. The first one
is the existence of global strong solutions to the NSC system
\eqref{nsft} for arbitrarily large initial data in $X$, provided that
the rotation speed $|\Omega|$ is sufficiently large (depending on the
data). To prove this, we closely follow the existence results that
have been established for rotating fluids in the whole space $\RR^3$,
see \cite[Chapter 5]{cdgg}. In particular, if the three-dimensional
part $\tilde u$ of the solution is not small at initial time, we assume
that the rotation speed $|\Omega|$ is large enough so that $\tilde u$
is rapidly damped by the dispersive effect of the linearized equation
\begin{equation}\label{freerot}
  \partial_t \tilde u +\Omega \PP(e_3 \wedge \tilde u) \,=\, 
  \Delta \tilde u~, \quad \div \tilde u \,=\, 0~.
\end{equation}
For the reader's convenience, we briefly recall in
Section~\ref{dispersec} and Appendix~\ref{dispersive} the Strichartz
estimates satisfied by the solutions of \eqref{freerot} with compact
support in Fourier space.  Except for the choice of the spatial
domain, the main difference of our approach with respect to
\cite{cdgg} is that we do not assume that the whole velocity field $u$
belongs to $L^2(\DD)^3$. As a consequence, we cannot use the energy
inequality which plays an important role in the classical approach. To
guarantee that the two-dimensional Navier-Stokes system has uniformly
bounded solutions, the hypothesis $\bar u_h=(\bar u_1,\bar u_2)^t \in 
L^2(\RR^2)^2$ is replaced by $\bar \omega_3 \in L^1(\RR^2)$, a
condition which allows for solutions with nonzero total circulation 
such as the Oseen vortex \eqref{Oseen1}, \eqref{Oseen2}.

The second part of Theorem~\ref{main}, which concerns the long-time
behavior of the solutions, is more in the spirit of the 
previous works \cite{thgcew2,roussier}. When stated more explicitly, 
our result shows that the solution $u(t,x,z)$ satisfies 
\[
  \|\tilde u(t)\|_{H^1(\DD)} + \|\bar u_3(t)\|_{H^1(\RR^2)}
  + \|\bar \omega_3(t)\|_{L^2(\RR^2)} \,\xrightarrow[t \to \infty]{}\, 
  0~,
\]
and
\begin{equation}
\label{omegaconv}
  \Bigl\|\bar \omega_3(t) - \frac{\alpha}{1+t}\,g\Bigl(\frac{\cdot}{
  \sqrt{1+t}}\Bigr)\Bigl\|_{L^1(\RR^2)} \,\xrightarrow[t \to \infty]{}\, 
  0~,
\end{equation}
where $g(\xi) = (4\pi)^{-1}\,e^{-|\xi|^2/4}$. In particular, if the 
total circulation $\alpha$ is nonzero, we see that $\bar \omega_3(t)$
does not converge to zero in the (scale invariant) space $L^1(\RR^2)$, 
but to the Oseen vortex with circulation $\alpha$, which is thus
the leading term in the asymptotic expansion of the solution as 
$t \to \infty$. This is in contrast with the case of finite-energy
solutions, which always converge to zero in the energy norm. 

\medskip
We conclude this introduction with a few additional remarks 
on the scope of Theorem~\ref{main}:

\noindent{\bf 1)} As is well-known, it is possible to prove the
existence of solutions to the NSC system \eqref{nsft} under weaker
assumptions on the initial data. For instance, it is sufficient 
to suppose that $\tilde u(0) \in H^{1/2}(\DD)^3$, $\bar u_3(0) 
\in L^2(\RR^2)$, and $\bar \omega_3(0) \in L^1(\RR^2)$, in which 
case the solution $u(t)$ will belong to $X$ for any positive 
time. Since we are mainly interested in the long-time behavior 
of the solutions, we disregard these technical details and prefer 
working directly in the (noncritical) space $X$.

\noindent{\bf 2)} Theorem~\ref{main} does not give any information on
the convergence rate towards Oseen's vortex. The proof shows that
$\|\nabla \bar u_3(t)\|_{L^2(\RR^2)} + \|\bar \omega_3(t)\|_{L^2(\RR^2)} 
= \cO(t^{-1/2})$ and $\|\tilde u(t)\|_{H^1(\DD)} = \cO(e^{-\nu t})$ 
for all $\nu < 4\pi^2$ as $t \to \infty$, but without additional 
assumptions on the data it is impossible to specify the decay rate 
of $\|\bar u_3(t)\|_{L^2(\RR^2)}$ or the convergence rate in
\eqref{omegaconv}. However, algebraic convergence rates can be
obtained if we assume that the initial data $\bar u_0(x)$ decay
sufficiently fast as $|x| \to \infty$, see \cite{thgcew2,roussier}.

\noindent{\bf 3)} In the proof of Theorem~\ref{main} we need a large
rotation speed $\Omega$ only to prove the existence of a global
solution, in the case where $\tilde u_0 = (1-Q)u_0$ is not small.
Once existence has been established, the convergence to Oseen's vortex
holds for any value of $\Omega$ and does not rely on the Coriolis
force at all. Since our domain $\DD$ has finite extension in the
vertical direction, we can use the Poincar\'e inequality to show that
$\tilde u(t)$ converges exponentially to zero as $t \to \infty$, but
this point is not crucial: Our proof can be adapted to cover the case
of the whole space $\RR^3$, if we assume as in \cite{cdgg} that $u =
\bar u + \tilde u$ with $\tilde u \in H^1(\RR^3)$, or even 
$\tilde u \in \dot H^{1/2}(\RR^3)$. In this situation the decay of 
$\tilde u(t)$ will not be exponential.

\noindent{\bf 4)} As is explained in \cite{roussier}, we can 
prove the analog of Theorem~\ref{main} in the layer 
$\RR^2 \times (0,1)$ with different bounday conditions, 
for instance stress-free conditions. The case of no-slip
(Dirichlet) boundary conditions is very different, because 
the solutions will converge exponentially to zero as 
$t \to \infty$, and the Oseen vortices can only appear as 
long-time transients.

\noindent{\bf 5)} A careful examination of the proof shows that 
the angular velocity $\Omega_0$ in Theorem~\ref{main} can be 
chosen in the following way:
\[
  \Omega_0 \,=\, \max\Bigl(K_0^2 \|\nabla \tilde u_0\|_{L^2} - K_0
  \,,\,0\Bigr)~, \quad \hbox{with}\quad K_0 \,=\, C e^{C\|u_0\|_X^8}~,
\]
where $\tilde u_0 = (1-Q)u_0$ and $C > 0$ is a universal constant. 
In particular, one can take $\Omega_0 = 0$ if $\tilde u_0$ is
sufficiently small, depending on $\bar u_0$. Of course, there
is no reason to believe that this result is sharp. 

\medskip The rest of this paper is organized as follows. In
Section~\ref{cauchy} we prove the existence part of Theorem~\ref{main}
using energy estimates for the full system \eqref{nsft} and dispersive
(Strichartz) estimates for the Rossby equation \eqref{freerot}.
Section~\ref{asymptotics} is devoted to the convergence proof, which
relies on a compactness argument and a transformation into
self-similar variables. In Appendix~\ref{BSlaw} we collect a few 
basic results concerning the Biot-Savart law in the domain $\DD$, 
and in Appendix~\ref{dispersive} we give a proof of the dispersive 
estimates for equation \eqref{freerot} which are used in the global 
existence proof.

\medskip\noindent
{\bf Acknowledgements.} The authors are indebted to Isabelle 
Gallagher for helpful discussions on several aspects of this work. 

\section{The Cauchy problem for the Navier-Stokes-Coriolis 
equation}
\label{cauchy}

In this section we prove that the Navier-Stokes-Coriolis system 
\eqref{nsft} is globally well-posed in the function space $X$ 
defined by \eqref{Xdef}, provided that the rotation speed 
$\Omega$ is sufficiently large depending on the initial data. 
The precise statement is:

\begin{thm}\label{glob}
For any initial data $u_0 \in X$, there exists $\Omega_0 \ge 0$ 
such that, for all $\Omega \in \RR$ with $|\Omega| \ge \Omega_0$, 
the NSC system \eqref{nsft} has a unique global solution
$u \in C^0([0,\infty),X)$ satisfying $u(0) = u_0$. Moreover, there 
exists $C > 0$ (depending on $u_0$) such that $\|u(t)\|_X \le C$ 
for all $t \ge 0$.
\end{thm}

As is clear from the proof, one can take $\Omega_0 = 0$ 
in Theorem~\ref{glob} (hence also in Theorem~\ref{main}) 
if the three-dimensional part $\tilde u_0 = (1-Q)u_0$ of 
the initial velocity field is sufficiently small in $X$, 
see Remark~\ref{smallrem} below. For large data, however, nobody 
knows how to prove global existence without assuming that 
the rotation speed $\Omega$ is large too. 

\subsection{Reformulation of the problem}
\label{equations}

If $u(t,x,z)$ is any solution of the NSC system \eqref{nsft}, we
decompose
\begin{equation}\label{dec1}
  u(t,x,z) \,=\, \bar u(t,x)+ \tilde u(t,x,z)~,
\end{equation}
where $\bar u = Qu$, $\tilde u = (1-Q)u$, and $Q$ is the vertical
average operator defined in \eqref{Qdef}. Our first task is to 
derive evolution equations for $\bar u$ and $\tilde u$. Integrating 
\eqref{nsft} over the vertical variable $z \in \TT^1$, and using
the fact that $\PP$ and $Q$ commute with each other (see 
Appendix~\ref{BSlaw}), we obtain
\begin{equation}\label{ns2d}
  \partial_t \bar{u}  + \PP[(\bar{u}\cdot\nabla)\bar{u} + Q
  (\tilde{u}\cdot\nabla)\tilde{u}] \,=\,\Delta \bar{u}~, \quad
  \div\bar{u} \,=\, 0~.
\end{equation}
This is a two-dimensional Navier-Stokes equation for the
three-component velocity field $\bar u(t,x)$, with a 
quadratic ``source term'' depending on $\tilde u$. Remark that 
the Coriolis force disappeared from \eqref{ns2d}, because $\curl 
(e_3 \wedge \bar u) = -\partial_z \bar u = 0$, so that 
$\PP(e_3 \wedge \bar u) = 0$. On the other hand, subtracting
\eqref{ns2d} from \eqref{nsft}, we find
\begin{equation}\label{ns3d}
  \partial_t \tilde{u} + \PP[(\bar{u}\cdot\nabla)\tilde{u} + 
  (\tilde{u}\cdot\nabla)\bar{u} + (1-Q)(\tilde{u}\cdot\nabla) 
  \tilde{u}] + \Omega \PP (e_3\wedge\tilde{u}) \,=\, 
  \Delta \tilde{u}~, \quad \div \tilde{u} \,=\, 0~.
\end{equation}
Thus $\tilde u(t,x,z)$ satisfies a three-dimensional 
Navier-Stokes-Coriolis system, which is linearly coupled to 
\eqref{ns2d} through the transport term $\PP(\bar{u}\cdot\nabla)
\tilde{u}$ and the stretching term $\PP(\tilde{u}\cdot\nabla)
\bar{u}$. 

As is explained in the introdution, the averaged velocity field 
$\bar{u}(t,x)$ can be represented by two scalar quantities, 
namely its vertical component $\bar{u}_3(t,x)$ and the third 
component $\bar{\omega}_3(t,x)$ of the averaged vorticity field.
Taking the third component of \eqref{ns2d} and using the 
fact that $(\PP\bar u)_3 = \bar u_3$ (see Appendix~\ref{BSlaw}), 
we obtain the following evolution equation: 
\begin{equation}\label{u3evol}
  \partial_t \bar{u}_3 + (\bar{u}_h \cdot \nabla)\bar{u}_3 + N_1 
  \,=\, \Delta \bar{u}_3~, \quad x \in \RR^2~, \quad t>0~,
\end{equation}
where $\bar u_h = (\bar u_1,\bar u_2)^t$ and $N_1 = Q(\tilde{u}\cdot
\nabla)\tilde{u}_3$. Similarly, if we take the third component
of \eqref{nsvor} and integrate the resulting equation over the 
vertical variable $z$, we find
\begin{equation}\label{omega3evol}
  \partial_t \bar{\omega}_3 + (\bar{u}_h \cdot \nabla)\bar{\omega}_3 
  + N_2 \,=\, \Delta\bar{\omega}_3~, \quad x \in \RR^2~, \quad t>0~,
\end{equation}
where $N_2 = Q((\tilde{u}\cdot\nabla)\tilde{\omega}_3 -
(\tilde{\omega}\cdot \nabla)\tilde{u}_3)$. Here we have used the
fact that $(\bar\omega \cdot\nabla)\bar u_3 = 0$, see 
\eqref{absys} below. 

By construction, the original NSC equation \eqref{nsft} is completely 
equivalent to the coupled system \eqref{ns3d}, \eqref{u3evol},
\eqref{omega3evol}. To prove local existence of solutions, we
consider the integral equations associated to these three 
PDE's (via Duhamel's formula), and we apply a standard fixed
point argument in the function space $C^0([0,T],X)$. 
The result is:

\begin{prop}\label{locex}
For any $r > 0$, there exists $T = T(r) > 0$ such that, 
for any $\Omega \in \RR$ and all initial data $u_0 \in X$
with $\|u_0\|_X \le r$, the Navier-Stokes-Coriolis system 
\eqref{nsft} has a unique local solution $u \in C^0([0,T],X)$ 
satisfying $u(0) = u_0$. 
\end{prop}

The proof of this statement uses classical arguments, which can be
found in \cite{fujitakato}, \cite{henry}, \cite{kato}, and will
therefore be omitted here. The fact that the local existence time $T$
depends on $u_0$ only through (an upper bound of) the norm $\|u_0\|_X$
is not surprising, because we work in a function space $X$ which is
not critical with respect to the scaling of the Navier-Stokes
equation. However, it is worth noticing that $T$ is independent of the
rotation speed $\Omega$. This is because the rotation does not act at
all on the two-dimensional part \eqref{u3evol}, \eqref{omega3evol} of
our system, whereas in \eqref{ns3d} it appears only in the term
$\Omega \PP(e_3 \wedge \tilde u)$, which is {\em skew-symmetric} in
the space $H^1(\DD)^3$ and therefore does not affect the estimates.

To prove global existence and conclude the proof of
Theorem~\ref{glob}, it remains to show that any solution $u \in
C^0([0,T],X)$ of \eqref{nsft} is bounded for all $t \in [0,T]$ by a
constant depending only on the initial data $u_0 = u(0)$. As is
well-known, this is relatively easy to do if the three-dimensional part
$\tilde u_0$ of the initial data is small in $H^1(\DD)$, see
\cite{fujitakato}, \cite{leray}. In the general case, we shall use the
dispersive properties of the Rossby equation \eqref{freerot} to prove
that the solution $\tilde u(t,x,z)$ of \eqref{ns3d} is rapidly damped
for positive times if the rotation speed $|\Omega|$ is sufficiently
large.

\subsection{Dispersive properties}
\label{dispersec}

Since our spatial domain $\DD = \RR^2 \times \TT^1$ is bounded in the
vertical direction, the Poincar\'e inequality implies that the
solutions of the linear equation \eqref{freerot} decay exponentially
to zero as $t \to \infty$. More precisely, for any $s \ge 0$ and all
divergence-free initial data $\tilde u_0 \in (1-Q)H^s(\DD)^3$, the
solution $\tilde u(t,x,z)$ of \eqref{freerot} satisfies
\begin{equation}\label{Hsest}
  \|\tilde u(t)\|_{H^s(\DD)} \,\le\, \|\tilde u_0\|_{H^s(\DD)}
  \,e^{-4\pi^2 t}~, \quad t \ge 0~.
\end{equation}
This estimate is straightforward to establish by computing the
time-derivative of $\|\tilde u(t)\|_{H^s}^2$ and using the Poincar\'e
inequality $\|\nabla \tilde u\|_{H^s}^2 \ge 4\pi^2 \|\tilde
u\|_{H^s}^2$ together with the fact that the Coriolis operator $\tilde
u \mapsto \PP(e_3 \wedge \tilde u)$ is skew-symmetric in $H^s(\DD)^3$
for divergence-free vector fields. Note in particular that
\eqref{Hsest} is independent ot the rotation speed $\Omega$.  However,
as is shown e.g. in \cite{cdgg}, additional information can be
obtained for large $|\Omega|$ if we exploit the dispersive effect of
the skew-symmetric term $\Omega \PP(e_3 \wedge \tilde u)$. The
corresponding Strichartz-type estimates are most conveniently derived
if we restrict ourselves to solutions with compact support in Fourier
space.

Throughout this paper, we use the following conventions for 
Fourier transforms. If $f \in L^2(\DD)$ or $L^2(\DD)^3$, we 
set 
\begin{equation}\label{fourier1}  
  f(x,z) \,=\, \frac{1}{2\pi}\int_{\RR^2} \sum_{n \in \ZZ} f_n(k)
  \,e^{i(k\cdot x +2\pi nz)} \dd k~,  \quad x \in \RR^2~, \quad z \in \TT^1~,
\end{equation}
where
\begin{equation}\label{fourier2}  
  f_n(k) \,=\, \frac{1}{2\pi} \int_{\RR^2}\int_{\TT^1} f(x,z)
  \,e^{-i(k\cdot x +2\pi nz)}\dd z \dd x~, \quad k \in \RR^2~, \quad
  n \in \ZZ~.  
\end{equation}
With these notations, the norm of $f$ in the Sobolev space 
$H^s(\DD)$ can be defined as
\begin{equation}\label{Hsnorm}
  \|f\|_{H^s} \,=\, \Bigl(\int_{\RR^2} \sum_{n \in \ZZ} 
  (1+|k|^2+4\pi^2 n^2)^s |f_n(k)|^2 \dd k\Bigr)^{1/2}~,
\end{equation}
where $|k|^2 = k_1^2 + k_2^2$. Given any $R > 0$, we denote by 
${\cal B}_R$ the ball
\begin{equation}\label{CRdef}
  {\cal B}_R \,=\, \Bigl\{(k,n) \in \RR^2 \times \ZZ \,\Big|\, 
  \sqrt{|k|^2+4\pi^2 n^2} \le R\Bigr\}~.
\end{equation}
Following closely the approach of \cite[Chap.~5]{cdgg}, we obtain
our main dispersion estimate: 

\begin{prop}\label{strichartz}
For any $R>0$, there exists $C_R>0$ such that, for all $\tilde u_0 
\in (1-Q)L^2(\DD)^3$ with $\div \tilde u_0 = 0$ and $\supp\,(\tilde 
u_0)_n(k) \subset {\cal B}_R$, the solution $\tilde u$ of 
\eqref{freerot} with initial data $\tilde u_0$ satisfies
\begin{equation}\label{strichest}
  \|\tilde u\|_{L^1(\RR_+,L^{\infty}(\DD))} \,\le\, C_R 
  |\Omega|^{-\frac{1}{4}} \|\tilde u_0\|_{L^2(\DD)}~.
\end{equation}
\end{prop}

For completeness, the proof of this proposition will be given in
Appendix~\ref{dispersive}. Estimate \eqref{strichest} clearly
demonstrates the dispersive effect of the Coriolis term in
\eqref{freerot}: If the initial data $\tilde u_0$ are compactly
supported in Fourier space, the $L^\infty$ norm of the solution
$\tilde u(t,\cdot)$ will be very small (for most values of time) if
the rotation speed $|\Omega|$ is large enough. This is in sharp contrast
with what happens for Sobolev norms, for which the best we can have is
estimate~\eqref{Hsest}. As a side remark, if we consider initial data
$\tilde u_0$ whose Fourier transform is supported {\em outside} the
ball ${\cal B}_R$, then we clearly have $\|\tilde u(t)\|_{H^s} \le 
\|\tilde u_0\|_{H^s}\,e^{-R^2 t}$ for all $t \ge 0$.

\medskip
Combining Proposition~\ref{strichartz} with estimate \eqref{Hsest}, 
we deduce the following useful corollary: 

\begin{cor}
\label{strichartz2}
Under the assumptions of Proposition~\ref{strichartz}, the solution 
$\tilde u$ of \eqref{freerot} satisfies, for any $p \in [1,+\infty]$ 
and any $q\in [2,+\infty]$ such that $\frac{1}{p}+\frac{2}{q} \le 1$,
\begin{equation}\label{strichest2}
  \|\tilde u\|_{L^p(\RR_+,L^q(\DD))} \,\le\, C_R 
  \langle\Omega\rangle^{-\frac{1}{4p}} \|\tilde u_0\|_{L^2(\DD)}~,
\end{equation}
where $\langle\Omega\rangle = (1 + |\Omega|^2)^{1/2}$. 
\end{cor}

\noindent\textbf{Proof.} Fix $s > 3/2$. Using Sobolev's embedding 
and our assumptions on $\tilde u_0$, we obtain from \eqref{Hsest}
\begin{equation}\label{cor1}
  \|\tilde u(t)\|_{L^\infty} \,\le\, C \|\tilde u(t)\|_{H^s}
  \,\le\, C\|\tilde u_0\|_{H^s}\,e^{-4\pi^2 t} \,\le\, 
  C_R \|\tilde u_0\|_{L^2}\,e^{-4\pi^2 t}~, \quad t \ge 0~,
\end{equation}
where $C_R$ denotes a generic positive constant depending only on $R$.
In particular, we have the estimate $\|\tilde u\|_{L^1(\RR_+,L^{\infty})}
\le C_R \|\tilde u_0\|_{L^2}$ for all $\Omega \in \RR$, so that
\eqref{strichest} holds with $|\Omega|$ replaced by $\langle\Omega
\rangle$. This gives \eqref{strichest2} for $(p,q) = (1,\infty)$, and
since the case $(p,q) = (\infty,\infty)$ is immediate from
\eqref{cor1}, we see that \eqref{strichest2} holds for all $p \in
[1,\infty]$ if $q = \infty$. Finally, as $\|\tilde u(t)\|_{L^2} \le
\|\tilde u_0\|_{L^2}$ for all $t \ge 0$, the general case follows by a
simple interpolation argument. \rule{2mm}{2mm}

\medskip To exploit the dispersive properties of the linear equation
\eqref{freerot} in the analysis of the nonlinear problem \eqref{ns3d},
we use the following decomposition, which is again borrowed from
\cite{cdgg}. Let $\chi \in C_0^\infty(\RR)$ be a cut-off function
satisfying $0 \le \chi(x) \le 1$ for all $x \in \RR$, $\chi(x) = 1$
for $|x| \le 1/2$ and $\chi(x) = 0$ for $|x| \ge 1$. Given any $R >
0$, we define the Fourier multiplyer ${\cal P}_R = \chi(|\nabla|/R)$
by the formula
\begin{equation}\label{PRdef}
  ({\cal P}_R f)_n(k) \,=\, \chi \Bigl(\frac{\sqrt{|k|^2 + 4\pi^2 
  n^2}}{R}\Bigr)f_n(k)~, \quad k \in \RR^2~, \quad n \in \ZZ~.
\end{equation}
If $\tilde u(t,x,z)$ is a solution of \eqref{ns3d} with initial 
data $\tilde u_0(x,z)$, we decompose
\begin{equation}\label{dec2}
  \tilde{u}(t,x,z) \,=\, \lambda(t,x,z) + r(t,x,z)~,
\end{equation}
where $\lambda(t,x,z)$ satisfies the linear Rossby equation
\begin{equation}\label{flre}
  \partial_t \lambda + \Omega \PP(e_3 \wedge \lambda) \,=\, 
  \Delta \lambda~, \quad \div \lambda \,=\, 0~,
\end{equation}
with initial data $\lambda_0 = {\cal P}_R \tilde{u}_0$. 
By construction, the remainder $r(t,x,z)$ is a solution 
of the nonlinear equation
\begin{equation}\label{req}
  \partial_t r +  \Omega \PP(e_3 \wedge r) + N_3 \,=\,\Delta r~, 
  \quad \div r \,=\, 0~,
\end{equation}
with initial data $r_0 = (1-{\cal P}_R)\tilde{u}_0$, where 
$N_3 = \PP [(\bar{u}\cdot\nabla)\tilde{u} + (\tilde{u}\cdot
\nabla)\bar{u} + (1-Q)(\tilde{u}\cdot\nabla)\tilde{u}]$. 

In the rest of this section, we consider equations~\eqref{flre},
\eqref{req} instead of \eqref{ns3d}, so that our final evolution
system consists of \eqref{u3evol}, \eqref{omega3evol}, \eqref{flre},
\eqref{req}. Given $u_0 = \bar u_0 + \tilde u_0\in X$, we will choose 
the parameter $R > 0$ large enough so that the initial data $r_0 =
(1-{\cal P}_R)\tilde{u}_0$ for equation~\eqref{req} are small in
$H^1(\DD)$. Then the rotation speed $|\Omega|$ will be taken large
enough so that we can exploit the dispersive estimates for
$\lambda(t,x,z)$ given by Corollary~\ref{strichartz2}.

\subsection{Energy estimates}
\label{energy}

We now derive the energy estimates which will be used to 
control the solutions of the nonlinear equations \eqref{u3evol}, 
\eqref{omega3evol}, \eqref{req}.  

\begin{prop}\label{energyprop}
There exist positive constants $C_0, C_1$ such that, if $u \in 
C^0([0,T],X)$ is a solution of \eqref{nsft} for some $\Omega \in 
\RR$, and if $u$ is decomposed as in \eqref{dec1}, \eqref{dec2}
for some $R > 0$, then the corresponding solutions of \eqref{u3evol}, 
\eqref{omega3evol}, \eqref{req} satisfy, for any $t \in (0,T]$:
\begin{align}
\label{sys1} 
  \frac{\dD}{\dD t} \|\bar{u}_3(t)\|_{L^2(\RR^2)}^2 \,&\le\, -\|\nabla 
  \bar{u}_3(t)\|_{L^2(\RR^2)}^2 +\|\tilde{u}(t)\|_{L^4(\DD)}^4~,\\
\label{sys2} 
  \frac{\dD}{\dD t} \|\nabla \bar{u}_3(t)\|_{L^2(\RR^2)}^2 \,&\le\, 
  -\|\Delta \bar{u}_3(t)\|_{L^2(\RR^2)}^2 + C_0(\|\nabla\bar{u}_3(t)
  \|_{L^2}^2\|\bar{\omega}_3(t)\|_{L^2}^2 + \| |\tilde{u}(t)|\,|\nabla
  \tilde{u}(t)|\|_{L^2}^2)~,\\
\label{sys3} 
  \frac{\dD}{\dD t} \|\bar{\omega}_3(t)\|_{L^2(\RR^2)}^2 \,&\le\, 
  -\|\nabla \bar{\omega}_3(t)\|_{L^2(\RR^2)}^2 + 8\| |\tilde{u}(t)| 
  |\nabla \tilde{u}(t)|\|_{L^2(\DD)}^2~, \\
\label{sys4} 
  \|\bar{\omega}_3(t)\|_{L^1(\RR^2)} \,&\le\, ~\|\bar{\omega}_3(0)
  \|_{L^1(\RR^2)} + 2\int_0^t \|\tilde{u}(s)\|_{L^2(\DD)}
  \|\Delta \tilde{u}(s)\|_{L^2(\DD)}\dd s~,\\
\label{sys5} 
  \frac{\dD}{\dD t} \|\nabla r(t)\|_{L^2(\DD)}^2 \,&\le\,  
  -\|\Delta r(t)\|_{L^2(\DD)}^2 + C_1 \|\nabla r(t)\|_{L^2}^2 \|\nabla 
  \bar{u}(t)\|_{L^2}^2 \|\Delta \bar{u}(t)\|_{L^2}^2  \\ \nonumber 
  \,&+\, C_1 (\|\bar{u}(t)\|_{L^4}^2 \|\nabla \lambda(t)\|_{L^4}^2 
  + \|\nabla \bar{u}(t)\|_{L^2}^2 \| \lambda(t)\|_{L^\infty}^2
  +\| |\tilde{u}(t)|\,|\nabla \tilde{u}(t)|\|_{L^2}^2)~.
\end{align}
\end{prop}

\begin{rem}\label{notations} Here and in what follows, if 
$f$ is a vector valued or matrix valued function, we denote by
$|f|$ the scalar function obtained by taking the Euclidean norm
of the entries of $f$. Given any $p \in [1,\infty]$, we define
$\|f\|_{L^p}$ as $\||f|\|_{L^p}$. With these conventions, if 
$\omega = \curl u$, we have for instance $|\omega| \le \sqrt{2}
|\nabla u|$ and $\|\omega\|_{L^2} = \|\nabla u\|_{L^2}$. 
\end{rem}

\medskip\noindent{\bf Proof.} To prove \eqref{sys1}, we multiply
both sides of \eqref{u3evol} by $\bar{u}_3$ and integrate 
over $\RR^2$. The transport term $(\bar u_h \cdot \nabla)\bar u_3$ 
gives no contribution, because $\bar u$ is divergence-free, and 
the diffusion term $\Delta \bar u_3$ produces the negative 
contribution $-\|\nabla \bar{u}_3\|_{L^2}^2$ after integrating 
by parts. Since
\[
  -\int_{\RR^2} \bar{u}_3 N_1 \dd x \,=\, -\int_{\DD} 
  \bar u_3 (\tilde u \cdot\nabla)\tilde u_3\dd x \dd z\,=\,  
  \int_{\DD} \tilde{u}_3(\tilde u \cdot \nabla \bar{u}_3)\dd x
  \dd z\,\le\, \frac12 \|\tilde u\|_{L^4}^4 + \frac12\|\nabla 
  \bar u_3 \|_{L^2}^2~,
\]
we obtain the desired estimate. In a similar way, to prove
\eqref{sys2}, we multiply \eqref{u3evol} by $-\Delta\bar{u}_3$ and
integrate over $\RR^2$. The transport term gives here a nontrivial
contribution which, after integrating by parts, can be bounded as
follows:
\begin{align*}
  \Bigl|\int_{\RR^2}(\Delta\bar{u}_3)(\bar{u}_h \cdot \nabla)
  \bar{u}_3 \dd x \Bigl| \,&\le\, \int_{\RR^2} |\nabla\bar{u}_3| 
  |\nabla \bar{u}_h| |\nabla \bar{u}_3| \dd x \,\le\, 
  \|\nabla \bar{u}_3\|_{L^4}^2 \|\nabla\bar{u}_h\|_{L^2} \\
  \,&\le\, C   \|\Delta\bar{u}_3\|_{L^2} \|\nabla \bar{u}_3\|_{L^2} 
  \|\bar{\omega}_3\|_{L^2} \,\le\, \frac{1}{4} \|\Delta \bar{u}_3
  \|_{L^2}^2 + C \|\nabla \bar{u}_3\|_{L^2}^2
  \|\bar{\omega}_3\|_{L^2}^2~.
\end{align*}
Here, to get from the first to the second line, we have used an 
interpolation inequality and the fact that $\bar u_h$ is 
obtained from $\bar\omega_3$ via the Biot-Savart law \eqref{BSbar2}, 
see Appendix~\ref{BSlaw}. Since we also have
\[
  \Bigl|\int_{\RR^2} \Delta \bar{u}_3 \,N_1 \dd x\Bigr| \,\le\, 
\int_{\DD} |\Delta \bar{u}_3| |\tilde u| |\nabla \tilde u|\dd x\dd z
  \,\le\, \frac{1}{4}\|\Delta \bar{u}_3\|_{L^2}^2 + \| |\tilde{u}| 
  |\nabla \tilde{u}|\|_{L^2}^2~,
\]
we obtain again the desired inequality.

On the other hand, multiplying \eqref{omega3evol} by $\bar{\omega}_3$ 
and integrating over $\RR^2$, we easily obtain \eqref{sys3}, because
\[
  -\int_{\RR^2} \bar{\omega}_3 N_2 \dd x \,=\, \int_{\DD}
  \Bigl(\tilde\omega_3 (\tilde u\cdot\nabla)\bar \omega_3 
  - \tilde u_3 (\tilde\omega \cdot\nabla)\bar \omega_3\Bigr)\dd x
  \dd z \,\le\, \frac12 \|\nabla\bar\omega_3\|_{L^2}^2 + 2 
  \||\tilde u| |\tilde \omega|\|_{L^2}^2~,
\]
and $|\tilde \omega|^2 \le 2|\nabla \tilde u|^2$. The prove 
\eqref{sys4} we observe that, since the vector field $\bar u_h$ 
is divergence-free, any solution of \eqref{omega3evol} in $L^1(\RR^2)$ 
satisfies
\[
  \|\bar{\omega}_3(t)\|_{L^1} \,\le\, \|\bar{\omega}_3(0)\|_{L^1}
  +\int_0^t \|N_2(s)\|_{L^1}\dd s~, \quad t \ge 0~.
\]
This bound can be established using the properties of the 
fundamental solution of the linear convection-diffusion 
equation $\partial_t f + (\bar u_h\cdot\nabla)f = \Delta f$, 
which will be recalled in Section~\ref{3.2} below. Since 
\[
  \|N_2\|_{L^1} \,\le\, \|\tilde u\|_{L^2}\|\nabla \tilde 
  \omega\|_{L^2} + \|\tilde\omega\|_{L^2}\|\nabla \tilde u
  \|_{L^2} \,\le\, \|\tilde u\|_{L^2}\|\Delta\tilde u\|_{L^2}
  + \|\nabla\tilde u\|_{L^2}^2 \,\le\, 2\|\tilde u\|_{L^2}\|
  \Delta\tilde u\|_{L^2}~,
\]
we obtain \eqref{sys4}. 

Finally, to prove \eqref{sys5}, we multiply \eqref{req} with
$-\Delta r$ and integrate over $\DD$. As was already explained, 
the Coriolis term $\Omega\PP(e_3\wedge r)$ gives no contribution, 
because it is skew-symmetric in any Sobolev space. So we just
have to bound the contributions of the nonlinear term $N_3$, which
are threefold. Since $\tilde u = \lambda + r$, the transport part
$\PP(\bar u\cdot\nabla)\tilde u$ in $N_3$ produces two terms, 
which can be estimated as follows:
\begin{align*}
  \Bigl| \int_{\DD} \Delta r & \cdot(\bar{u} \cdot \nabla)\lambda\dd x \dd z 
  \Bigr| \,\le\, \frac{1}{10} \|\Delta r\|_{L^2}^2 + C 
  \|\bar{u}\|_{L^4}^2 \|\nabla \lambda\|_{L^4}^2~,\\
  \Bigl| \int_{\DD} \Delta r & \cdot(\bar{u} \cdot \nabla) r\dd x \dd z\Bigr|
  \,\le\, \int_{\DD} |\nabla r| |\nabla \bar{u}| |\nabla r|\dd x \dd z
  \,\le\, \|\nabla r\|_{L^{\frac{8}{3}}}^2 \|\nabla \bar{u}\|_{L^4}\\
  \,&\le\, C\|\nabla r\|_{L^2}^{\frac{5}{4}} \|\Delta r\|_{L^2}^{\frac{3}{4}}
  \|\nabla \bar{u}\|_{L^2}^{\frac{1}{2}} \|\Delta \bar{u}\|_{L^2}^{\frac{1}{2}}
  \,\le\, C\|\nabla r\|_{L^2}^{\frac{1}{2}}\|\Delta r\|_{L^2}^{\frac{3}{2}}  
  \|\nabla \bar{u}\|_{L^2}^{\frac{1}{2}} \|\Delta \bar{u}\|_{L^2}^{\frac{1}{2}}\\
  \,&\le\, \frac{1}{10} \|\Delta r\|_{L^2}^2 + C \|\nabla r\|_{L^2}^2  
  \|\nabla \bar{u}\|_{L^2}^2 \|\Delta \bar{u}\|_{L^2}^2~.
\end{align*} 
Here we have used interpolation inequalities, Sobolev embeddings, 
and the Poincar\'e inequality $\|\nabla r\|_{L^2} \le C\|\Delta
r\|_{L^2}$. The two terms produced by the stretching part $\PP(\tilde 
u\cdot\nabla)\bar u$ in $N_3$ can be estimated in a similar way:
\begin{align*}
  \Bigl|\int_{\DD} \Delta r \cdot(\lambda\cdot \nabla)\bar{u}\dd x\dd z\Bigr| 
  \,&\le\, \frac{1}{10} \|\Delta r\|_{L^2}^2 +C \|\lambda\|_{L^{\infty}}^2
  \|\nabla \bar{u}\|_{L^2}^2~,\\
  \Bigl|\int_{\DD} \Delta r \cdot(r \cdot \nabla) \bar{u}\dd x \dd z\Bigr|
  \,&\le\, \|\Delta r\|_{L^2} \|r\|_{L^4} \|\nabla \bar{u}\|_{L^4}
  \,\le\, C\|\Delta r\|_{L^2} \|\nabla r\|_{L^2} 
  \|\nabla \bar{u}\|_{L^2}^{\frac{1}{2}} \|\Delta \bar{u}\|_{L^2}^{\frac{1}{2}} \\
  \,&\le\, \frac{1}{10} \|\Delta r\|_{L^2}^2 + C \|\nabla r\|_{L^2}^2 
\|\nabla \bar{u}\|_{L^2}^2 \|\Delta \bar{u}\|_{L^2}^2~. 
\end{align*}
Finally, the contribution of the quadratic term $\PP(1-Q)(\tilde
u\cdot\nabla)\tilde{u}$ in $N_3$ satisfies
\[
 \Bigl|\int_{\DD} \Delta r \cdot(\tilde u\cdot \nabla)\tilde{u} \dd x\dd z
 \Bigr| \,\le\, \frac{1}{10}\|\Delta r\|_{L^2}^2 + C \| |\tilde{u}
  |\,|\nabla \tilde{u}|\|_{L^2}^2~.
\]
Collecting all these estimates, we obtain \eqref{sys5}. This 
concludes the proof. \rule{2mm}{2mm}

\subsection{Global existence}
\label{global}

In this section, we combine the dispersive properties of 
Section~\ref{dispersec} and the energy estimates of 
Section~\ref{energy} to complete the proof of Theorem~\ref{glob}.
We start with a preliminary result, which summarizes in a 
convenient way four of the five inequalities established 
in Proposition~\ref{energyprop}. 

\begin{lm}\label{intermediaire}
There exist positive constants $C_2$, $C_3$, and $C_4$ such 
that the following holds. Let $u \in C^0([0,T],X)$ be a solution of
\eqref{nsft} for some $\Omega \in \RR$, which is decomposed as in 
\eqref{dec1}, \eqref{dec2} for some $R > 0$. Assume moreover
that there exist $K \ge 1$ and $\epsilon \in (0,1]$ such that
the corresponding solutions of \eqref{u3evol}, \eqref{omega3evol}, 
\eqref{req} satisfy
\begin{equation}\label{apriori}
  \|\nabla \bar{u}_3(t)\|_{L^2(\RR^2)} \,\le\, K^2~, \quad
  \|\bar{\omega}_3(t)\|_{L^2(\RR^2)} \,\le\, K~, \quad
  \|\nabla r(t)\|_{L^2(\DD)} \,\le\, \epsilon~,
\end{equation}
for all $t \in [0,T]$. If we define
\begin{equation}\label{Phidef}
  \Phi(t) \,=\, \|\bar{u}_3(t)\|_{L^2(\RR^2)}^2 + \|\bar{\omega}_3(t)
  \|_{L^2(\RR^2)}^2 +\delta \|\nabla \bar{u}_3\|_{L^2(\RR^2)}^2+ 
  \|\nabla r(t)\|_{L^2(\DD)}^2~,
\end{equation}
for some $\delta \in (0,1]$, then 
\begin{align}\nonumber
  \frac{\dD}{\dD t}\Phi(t) \,\le\, 
  &-(\|\nabla \bar{u}_3(t)\|_{L^2}^2 + \|\nabla \bar{\omega}_3(t)
  \|_{L^2}^2 + \delta \|\Delta \bar{u}_3(t)\|_{L^2}^2 + 
  \|\Delta r(t)\|_{L^2}^2) \\ \label{Phidot}
  \,&+\, C_0\delta K^2 \|\nabla \bar{u}_3(t)\|_{L^2}^2
  + C_2 \epsilon^2 K^4 \|\Delta \bar{u}(t)\|_{L^2}^2 
  + C_3 \epsilon^2 \|\Delta r(t)\|_{L^2}^2 \\ \nonumber
  \,&+\, (\delta^{-1}\Phi(t) + K\|\bar{\omega}_3(t)\|_{L^1})G(t)
  + F(t) + \epsilon^2 G(t)~,
\end{align}
for all $t \in (0,T]$, where 
\begin{align}\label{FGdef}
  F(t) \,&=\, C_4(\|\lambda(t)\|_{L^4}^4 + \|\lambda(t)\|_{L^\infty}^2
  \|\nabla\lambda(t)\|_{L^2}^2)~, \\ \nonumber
  G(t) \,&=\, C_4(\|\nabla\lambda(t)\|_{L^\infty}^2 + \|\lambda(t)
  \|_{L^\infty}^2 + \|\nabla\lambda(t)\|_{L^4}^2)~.
\end{align}
\end{lm}

\medskip\noindent{\bf Proof.} If $\Phi$ is defined by 
\eqref{Phidef}, it follows immediately from 
Proposition~\ref{energyprop} that
\begin{align}\nonumber
\frac{\dD}{\dD t} \Phi(t)  \,&\le\, 
  -\, (\|\nabla \bar{u}_3(t)\|_{L^2}^2 +\|\nabla \bar{\omega}_3(t)
  \|_{L^2}^2 + \delta \|\Delta \bar{u}_3(t)\|_{L^2}^2 + \|\Delta 
  r(t)\|_{L^2}^2) \\ \label{Phidot2}
  \,&+\, \|\tilde{u}(t)\|_{L^4}^4 + C \||\tilde{u}(t)| |\nabla 
  \tilde{u}(t)|\|_{L^2}^2 + C_0 \delta \|\nabla \bar{u}_3(t)
  \|_{L^2}^2 \|\bar{\omega}_3(t)\|_{L^2}^2\\ \nonumber
  \,&+\, C (\|\nabla r(t)\|_{L^2}^2 \|\nabla \bar{u}(t)\|_{L^2}^2 
  \|\Delta \bar{u}(t)\|_{L^2}^2 
  + \|\bar{u}(t)\|_{L^4}^2 \|\nabla \lambda (t)\|_{L^4}^2 
  + \|\nabla \bar{u}(t)\|_{L^2}^2 \|\lambda(t)\|_{L^{\infty}}^2)~.
\end{align}
Using interpolation inequalities, Sobolev embeddings, and the a priori
bounds \eqref{apriori}, we first get
\begin{align*}
  \|\tilde{u}(t)\|_{L^4}^4  \,&\le\, C(\|r(t)\|_{L^4}^4 +
  \|\lambda(t)\|_{L^4}^4) \,\le\, C(\|r(t)\|_{L^6}^3 \|r(t)\|_{L^2} + 
  \|\lambda(t)\|_{L^4}^4)\\
  \,&\le\, C(\|\nabla r(t)\|_{L^2}^3 \|r(t)\|_{L^2} + \|\lambda(t)
  \|_{L^4}^4) \,\le\, C\epsilon^2\|\nabla r(t)\|_{L^2} \|r(t)\|_{L^2} 
  + F_1(t)~,
\end{align*}
where $F_1(t) = C\|\lambda(t)\|_{L^4}^4$. Proceeding in the same way, 
we also obtain
\begin{align*}
  \||\tilde{u}||\nabla \tilde{u}|\|_{L^2}^2 \,&\le\, 
  C(\||r||\nabla r|\|_{L^2}^2 + \||r||\nabla\lambda|\|_{L^2}^2 
  + \||\lambda||\nabla r|\|_{L^2}^2 + \||\lambda||\nabla
  \lambda|\|_{L^2}^2)\\
  \,&\le\, C (\|r\|_{L^6}^2 \|\nabla r \|_{L^3}^2 +
  \|r\|_{L^2}^2\|\nabla\lambda\|_{L^\infty}^2 + 
  \|\nabla r\|_{L^2}^2\|\lambda\|_{L^\infty}^2 + 
  \|\lambda\|_{L^\infty}^2\|\nabla\lambda\|_{L^2}^2)~,
\end{align*}
so that $\||\tilde{u}(t)| |\nabla \tilde{u}(t)|\|_{L^2}^2 \le 
C\epsilon^2\|\nabla r(t)\|_{L^2}\|\Delta r(t)\|_{L^2} + \epsilon^2
G_1(t) + F_2(t)$, 
where 
\[
  G_1(t) \,=\, C(\|\nabla\lambda(t)\|_{L^{\infty}}^2 +
  \|\lambda(t)\|_{L^{\infty}}^2)~, \quad F_2(t) \,=\, 
  C\|\lambda(t)\|_{L^\infty}^2 \|\nabla \lambda(t)\|_{L^2}^2~.
\]
It remains to estimate the last four terms in the right-hand side
of \eqref{Phidot2}. The first two in this group are independent 
of $\lambda$, and are simply bounded using assumption \eqref{apriori} 
and the fact that $\|\nabla \bar u\|_{L^2}^2 = \|\nabla \bar 
u_3\|_{L^2}^2 + \|\bar \omega_3\|_{L^2}^2$. On the other hand, in view 
of Proposition~\ref{BSbarprop}, we have
\[
  \|\bar{u}\|_{L^4}^2 \,\le\, C(\|\bar{u}_3\|_{L^4}^2 + 
  \|\bar{\omega}_3\|_{L^{\frac{4}{3}}}^2) \,\le\, 
  C(\|\bar{u}_3\|_{L^2}\|\nabla \bar{u}_3\|_{L^2} + 
  \|\bar{\omega}_3\|_{L^1} \|\bar{\omega}_3\|_{L^2})~,   
\]
hence 
\[
  \|\bar{u}(t)\|_{L^4}^2 \|\nabla \lambda (t)\|_{L^4}^2
  \,\le\, (\delta^{-1/2}\Phi(t) + K\|\bar{\omega}_3(t)\|_{L^1})
  G_2(t)~,
\]
where $G_2(t) = C\|\nabla\lambda (t)\|_{L^4}^2$. Similarly, we find
$\|\nabla\bar{u}(t)\|_{L^2}^2 \|\lambda(t) \|_{L^{\infty}}^2 \le
\delta^{-1}\Phi(t)G_1(t)$. Thus, using the Poincar\'e inequality
$\|r\|_{L^2} \le \|\nabla r\|_{L^2} \le \|\Delta r\|_{L^2}$, 
we see that \eqref{Phidot} holds with $F(t) = F_1(t) + F_2(t)$ 
and $G(t) = G_1(t) + G_2(t)$.  \rule{2mm}{2mm}

\begin{rem}\label{FGrem} In view of Corollary~\ref{strichartz2}, 
there exists a constant $C_R > 0$ (depending only on $R$) such that
\begin{equation}\label{FGbdd}
  \int_0^\infty F(t)\dd t \,\le\, C_R\langle\Omega\rangle^{-\frac{1}{4}}
  \|\tilde{u}_0\|_{L^2}^4~, \quad \hbox{and}\quad
  \int_0^\infty G(t)\dd t \,\le\, C_R\langle\Omega\rangle^{-\frac{1}{4}}
  \|\tilde{u}_0\|_{L^2}^2~.
\end{equation}
\end{rem}

\begin{rem}\label{Constrem}
Without loss of generality, we shall assume henceforth that 
the constants which appear in Proposition~\ref{energyprop} 
and Lemma~\ref{intermediaire} satisfy $C_i \ge 1$, $i = 0,\dots,4$. 
\end{rem}

\medskip\noindent
{\bf Proof of theorem \ref{glob}.} Given $u_0 \in X$, we define 
$\bar u(0) = Qu_0$, $\tilde u_0 = (1-Q)u_0$, and $\bar \omega(0) = 
\curl(Qu_0)$, where $Q$ is the vertical average operator 
\eqref{Qdef}. We first choose $K \ge 1$ such that
\begin{equation}\label{Kchoice}
  \|\bar u_3(0)\|_{H^1(\RR^2)}^2 + \|\bar\omega_3(0)\|_{L^2(\RR^2)}^2
  + \|\bar\omega_3(0)\|_{L^1(\RR^2)} + 2\|\tilde u_0\|_{H^1(\DD)}^2 
  \,\le\, \frac{K^2}{16 C_0}~,
\end{equation}
where $C_0 \ge 1$ is as in Proposition~\ref{energyprop}. 
Next, we take $\epsilon \in (0,1]$ sufficiently small so that
\begin{equation}\label{epsilonchoice}
  \epsilon^2 \,\le\, \min\Bigl\{\frac{1}{2C_3}\,,\,
  \frac{\delta}{2C_2 K^4}\Bigr\}~, \quad \hbox{where} \quad
  \delta \,=\, \frac{1}{2C_0 K^2} \,\in\, (0,1]~, 
\end{equation}
and $C_2 \ge 1$, $C_3 \ge 1$ are as in Lemma~\ref{intermediaire}.
Once this is done, we set $\lambda_0 = {\cal P}_R \tilde u_0$ and
$r_0 = (1-{\cal P}_R) \tilde u_0$, where ${\cal P}_R$ is the 
Fourier localization operator defined by \eqref{PRdef}. We assume 
that the parameter $R > 0$ is sufficiently large so that 
\begin{equation}\label{Rchoice}
  4\,e^{2C_1 K^8}\|\nabla r_0\|_{L^2}^2 \,\le\, \epsilon^2~,
\end{equation}
and we denote by $\lambda(t,x,z)$ the solution of \eqref{flre} 
with initial data $\lambda_0$. Finally, using Remark~\ref{FGrem}, 
we choose $\Omega_0 \ge 0$ sufficiently large so that, if $|\Omega| 
\ge \Omega_0$,
\begin{equation}\label{Omegachoice1}
  \int_0^\infty G(t)\dd t \,\le\, \delta\log(2)~, \quad
  \int_0^\infty (F(t) + \epsilon^2 G(t))\dd t \,\le\, 
  \frac{K^2}{16 C_0}~,
\end{equation}
and 
\begin{equation}\label{Omegachoice2}
  4\,e^{2C_1 K^8}\int_0^\infty \Bigl(F(t) + (K^4+\epsilon^2)G(t)
  \Bigr)\dd t \,\le\, \epsilon^2~.
\end{equation}

\begin{rem}\label{smallrem}
If $\tilde u_0$ is small enough so that $4\,e^{C_1 K^8}\|\nabla \tilde
u_0\|_{L^2}^2 \le \epsilon^2$, then we can take formally $R = 0$, 
so that $r_0 = \tilde u_0$ and $\lambda_0 = 0$. In that case, 
one has $F(t) = G(t) \equiv 0$, and \eqref{Omegachoice1}, 
\eqref{Omegachoice2} are of course satisfied for any $\Omega 
\in \RR$.
\end{rem}

By Proposition~\ref{locex}, equation \eqref{nsft} has a unique maximal 
solution $u \in C^0([0,T_*),X)$ with initial data $u_0$, where 
$T_* \in (0,+\infty]$ denotes the maximal existence time. If we 
decompose $u(t) = \bar u(t) + \lambda(t) + r(t)$ as in \eqref{dec1}, 
\eqref{dec2}, then $\bar u_3(t)$, $\bar \omega_3(t)$, $r(t)$ 
are solutions of \eqref{u3evol}, \eqref{omega3evol}, \eqref{req}, 
respectively, and we know from \eqref{Kchoice} and \eqref{Rchoice} 
that 
\[
  \|\nabla \bar{u}_3(0)\|_{L^2(\RR^2)} \,\le\, \frac{K}{4}~, \quad
  \|\bar{\omega}_3(0)\|_{L^2(\RR^2)} \,\le\, \frac{K}{4}~, \quad
  \|\nabla r_0\|_{L^2(\DD)} \,\le\, \frac{\epsilon}{2}~.
\]
Thus, by continuity, the bounds \eqref{apriori} will be satisfied 
at least for $t > 0$ sufficiently small. Let
\begin{equation}\label{Tdef}
  T \,=\, \sup\Bigl\{\tilde{T} \in [0,T_*) \,\Big|\, \hbox{The 
  bounds }\eqref{apriori}\hbox{ hold for all }t\in [0,\tilde{T}]
  \Bigr\} \,\in\, (0,T_*]~.
\end{equation}
We shall prove that $T = T_*$. This implies of course that $T = 
T_* = +\infty$, and that the solution $u(t)$ of \eqref{nsft} 
stays bounded in $X$ for all $t \ge 0$, as is claimed in 
Theorem~\ref{glob}. 

Assume on the contrary that $0 < T < T_*$, and let $\Psi(t)
= \Phi(t) + \|\bar{\omega}_3(t)\|_{L^1}$, where $\Phi$ is 
defined in \eqref{Phidef}. Using \eqref{Phidot} and 
\eqref{epsilonchoice}, we find
\begin{align}\label{last1}
  \Phi(t) \,&+\, \frac12 \int_0^t \Bigl(\|\nabla \bar{u}_3(s)\|_{L^2}^2 
  + \|\nabla \bar{\omega}_3(s)\|_{L^2}^2 + \delta\|\Delta\bar{u}_3(s)
  \|_{L^2}^2 + \|\Delta r(s)\|_{L^2}^2\Bigr)\dd s \\ \nonumber
  \,&\le\, \Phi(0) + \delta^{-1}\int_0^t \Psi(s)G(s)\dd s + 
  \int_0^t (F(s)+\epsilon^2 G(s))\dd s~, \quad t \in [0,T]~.
\end{align}
On the other hand, since
\[
  2\|\tilde{u}(t)\|_{L^2}\|\Delta\tilde{u}(t)\|_{L^2} \,\le\, 
  \frac{1}{2\pi^2}\|\Delta\tilde{u}(t)\|_{L^2}^2 \,\le\, 
  \frac{1}{\pi^2}\Bigl(\|\Delta r(t)\|_{L^2}^2 + 
  \|\Delta\lambda(t)\|_{L^2}^2\Bigr)~,
\]
it follows from \eqref{sys4} that
\begin{equation}\label{last2}
  \|\bar{\omega}_3(t)\|_{L^1} \,\le\, \|\bar{\omega}_3(0)\|_{L^1} 
  + \frac14 \int_0^t \|\Delta r(s)\|_{L^2}^2\dd s + 
  \|\nabla\tilde u_0\|_{L^2}^2~, \quad t \in [0,T]~. 
\end{equation}
Here we have used the fact that $2\int_0^\infty \|\Delta 
\lambda(t)\|_{L^2}^2\dd t = \|\nabla\lambda_0\|_{L^2}^2 \le 
\|\nabla\tilde u_0\|_{L^2}^2$ by \eqref{flre}. Summing up
\eqref{last1} and \eqref{last2}, 
we obtain for $t \in [0,T]$:
\begin{align}\label{last3}
  \Psi(t) \,&+\, \frac12 \int_0^t \Bigl(\|\nabla \bar{u}_3(s)\|_{L^2}^2 
  + \|\nabla \bar{\omega}_3(s)\|_{L^2}^2 + \delta\|\Delta\bar{u}_3(s)
  \|_{L^2}^2 + \frac12\|\Delta r(s)\|_{L^2}^2\Bigr)\dd s \\ \nonumber
  \,&\le\, \Psi(0) + \|\nabla\tilde u_0\|_{L^2}^2 +
  \delta^{-1}\int_0^t \Psi(s)G(s)\dd s + \int_0^t (F(s)+\epsilon^2 
  G(s))\dd s~.
\end{align}
This integral inequality for $\Psi(t)$ can be integrated 
using Gronwall's lemma. In view of \eqref{Kchoice}, \eqref{Rchoice}
and \eqref{Omegachoice1}, we easily obtain
\begin{align}\label{last4}
  \Psi(t) \,&+\, \frac12 \int_0^t \Bigl(\|\nabla \bar{u}_3(s)\|_{L^2}^2 
  + \|\nabla \bar{\omega}_3(s)\|_{L^2}^2 + \delta\|\Delta\bar{u}_3(s)
  \|_{L^2}^2 + \frac12\|\Delta r(s)\|_{L^2}^2\Bigr)\dd s \\ \nonumber
  \,&\le\, 2\Bigl(\Psi(0) + \|\nabla\tilde u_0\|_{L^2}^2 + 
  \int_0^t (F(s)+\epsilon^2 G(s))\dd s\Bigr) \,\le\, 
  \frac{K^2}{4C_0}~,
\end{align}
for all $t \in [0,T]$. In a similar way, using  \eqref{sys5}, 
\eqref{apriori} and proceeding as in the proof of 
Lemma~\ref{intermediaire}, we find
\begin{align}\label{last5}
  \|\nabla r(t)\|_{L^2}^2 \,&+\, \frac12\int_0^t \|\Delta r(s)\|_{L^2}^2 
  \dd s \,\le\, \|\nabla r_0\|_{L^2}^2 + 2 C_1 K^4 \int_0^t 
  \|\nabla r(s)\|_{L^2}^2 \|\Delta \bar{u}(s)\|_{L^2}^2 \dd s \\ \nonumber
  \,&+\, \delta^{-1}\int_0^t \Psi(s)G(s)\dd s + 
    \int_0^t (F(s)+\epsilon^2 G(s))\dd s~, \quad t \in [0,T]~.
\end{align}
From \eqref{last4} we know that $\int_0^t \|\Delta\bar{u}(s)\|_{L^2}^2
\dd s \le 2K^2/(4C_0\delta) = K^4$. Thus we can apply Gronwall's
lemma to \eqref{last5} and, using in addition \eqref{Rchoice} and
\eqref{Omegachoice2}, we obtain
\begin{align}\label{last6}
  \|\nabla r(t)\|_{L^2}^2 \,&+\, \frac12\int_0^t \|\Delta r(s)\|_{L^2}^2 
  \dd s \\ \nonumber 
  \,&\le\, e^{2C_1 K^8}\Bigl(\|\nabla r_0\|_{L^2}^2 + K^4 
  \int_0^t G(s)\dd s + \int_0^t (F(s)+\epsilon^2 G(s))
  \dd s\Bigr) \,\le\, \frac{\epsilon^2}{2}~,
\end{align}
for all $t \in [0,T]$. 

Now, it follows immediately from \eqref{last4}, \eqref{last6}
that
\[
  \|\nabla\bar{u}_3(t)\|_{L^2}^2 \,\le\, \frac{K^4}{2}~, \quad
  \|\bar{\omega}_3(t)\|_{L^2}^2 \,\le\, \frac{K^2}{4C_0}~, \quad
  \|\nabla r(t)\|_{L^2}^2 \,\le\, \frac{\epsilon^2}{2}~, 
\]
for all $t \in [0,T]$, which obviously contradicts the definition 
\eqref{Tdef} of $T$. Thus $T = T_* = +\infty$, and estimates 
\eqref{apriori}, \eqref{last4}, \eqref{last6} hold for all
$t \ge 0$. This concludes the proof of Theorem~\ref{glob}. 
\rule{2mm}{2mm}

\section{Convergence to Oseen Vortices}
\label{asymptotics}

To complete the proof of Theorem~\ref{main}, it remains to show that
the global solution $u(t,x,z)$ of the Navier-Stokes-Coriolis system
\eqref{nsft} constructed in Section~\ref{cauchy} converges to Oseen's
vortex as $t \to \infty$. To do that, we decompose $u(t,x,z) = \bar
u(t,x) + \tilde u(t,x,z)$ as in \eqref{dec1}, and we first show that
the three-dimensional part $\tilde u(t)$ converges exponentially to zero
in $H^1(\DD)^3$, due to Poincar\'e's inequality. We next turn our
attention to the two-dimensional part $\bar u$, and prove that the
third component $\bar u_3(t)$ decays to zero in $H^1(\RR^2)$.  Finally,
the most delicate point is to show that $\bar \omega_3(t)$ converges
to Oseen's vortex in $L^1(\RR^2)$ as $t \to \infty$.  Here the main
ingredients are a transformation into self-similar variables, a
compactness estimate for the rescaled solution, and a characterization
of the complete trajectories of the two-dimensional Navier-Stokes
equation which was obtained in \cite{thgcew2}.

\subsection{Exponential decay of $\tilde u$}
\label{3.1}

We recall from \eqref{dec2} that $\tilde u(t,x,z) = r(t,x,z) +
\lambda(t,x,z)$, where $\lambda$ satisfies the linear equation
\eqref{flre} and $r$ is a solution of \eqref{req}. We already know 
that $\|\lambda(t)\|_{H^s} \le C\,e^{-4\pi^2 t}$ for all $t \ge 0$
and any $s \ge 0$, see \eqref{Hsest}, so it remains to estimate
$r(t,x,z)$. We start from equation \eqref{sys5} which, in view of the
global bound obtained in Theorem~\ref{glob} and the estimate above for
$\lambda$, implies
\begin{equation}\label{rineq}
  \frac{\dD}{\dD t}\|\nabla r(t)\|_{L^2}^2 + \frac12 \|\Delta r(t)
  \|_{L^2}^2 \,\le\, C_1 \|\nabla r(t)\|_{L^2}^2 \|\Delta \bar u(t)
  \|_{L^2}^2 + C_2\,e^{-8\pi^2 t}~,
\end{equation}
for some constants $C_1, C_2 > 0$ (depending on the initial data). 
Fix $0 < \mu \le 2\pi^2$ and let $f(t) = e^{\mu t} \|\nabla
r(t)\|_{L^2}^2$.  Using \eqref{rineq} and the Poincar\'e inequality
$\|\Delta r\|_{L^2} \ge 2\pi \|\nabla r\|_{L^2}$, we find
\begin{align}\nonumber
  f'(t) \,&\le\, e^{\mu t}\Bigl(\mu\|\nabla r(t)\|_{L^2}^2 
  -\frac12 \|\Delta r(t)\|_{L^2}^2 + C_1 \|\nabla r(t)\|_{L^2}^2 
  \|\Delta \bar u(t)\|_{L^2}^2 + C_2\,e^{-8\pi^2 t}\Bigr) \\
  \label{fineq}
  \,&\le\, C_1 f(t)\|\Delta \bar u(t)\|_{L^2}^2 + C_2\,e^{-(8\pi^2
   -\mu) t}~.
\end{align}
Since $\int_0^\infty \|\Delta \bar u(t)\|_{L^2}^2\dd t < \infty $
by \eqref{last4}, it follows from \eqref{fineq} that $f(t) \le C_3$ 
for all $t \ge 0$, hence $\|\nabla r(t)\|_{L^2} \le C_3\,e^{-\mu t/2}$ 
for some $C_3 > 0$. As $\|\tilde u\|_{H^1} \approx \|\nabla \tilde u
\|_{L^2} \le \|\nabla r\|_{L^2} + \|\nabla \lambda\|_{L^2}$, 
this proves that $\tilde u(t)$ converges exponentially to zero 
in $H^1(\DD)^3$ as $t \to \infty$. The decay rate we have obtained 
so far is not optimal, but it is sufficient to conclude the 
proof of Theorem~\ref{main}. 

To get the optimal decay rate, the simplest solution is to go back 
to equation \eqref{ns3d} satisfied by $\tilde u$. Using straightforward
estimates to bound the nonlinear terms, we arrive at the differential
inequality
\begin{align}\nonumber
  \frac{\dD}{\dD t}\|\nabla \tilde u(t)\|_{L^2}^2 \,&\le\, 
  -2 \|\Delta\tilde u(t)\|_{L^2}^2 + \int_\DD \Delta\tilde u(t)\cdot
  N_3(t) \dd x\dd z \\ \label{uineq}
  \,&\le\, -2 \|\Delta\tilde u(t)\|_{L^2}^2 + C\|\Delta\tilde u(t)\|_{L^2}
  \|\nabla\tilde u(t)\|_{L^2}(\|\nabla\tilde u(t)\|_{L^3} + 
  \|\nabla\bar u(t)\|_{L^3})~,
\end{align}
where $C > 0$ is a universal constant. Now we observe that 
\begin{equation}\label{finint}
  \int_0^\infty (\|\nabla\tilde u(t)\|_{L^3}^2 + \|\nabla\bar 
  u(t)\|_{L^3}^2)\dd t \,\le\, C \!\int_0^\infty (\|\nabla\tilde 
  u(t)\|_{L^3}^2 + \|\nabla\bar u_3(t)\|_{L^3}^2 + 
  \|\bar \omega_3(t)\|_{L^3}^2)\dd t \,<\, \infty~.
\end{equation}
For $\tilde u$ and $\bar u_3$, this claim follows \eqref{last4}, 
\eqref{last5}, because $\|\nabla\tilde u\|_{L^3}^2 \le C\|\Delta 
\tilde u\|_{L^2}^2 \le C(\|\Delta r\|_{L^2}^2 + \|\Delta 
\lambda\|_{L^2}^2)$, and $\|\nabla\bar u_3\|_{L^3}^2 \le 
C\|\nabla\bar u_3\|_{L^2}^{4/3}\|\Delta \bar u_3\|_{L^2}^{2/3} \le 
C(\|\nabla\bar u_3\|_{L^2}^2 + \|\Delta \bar u_3\|_{L^2}^2)$.
On the other hand, the decay rates established in 
Section~\ref{3.3} below will show that $\|\bar \omega_3(t)\|_{L^3} =
\cO(t^{-2/3})$ as $t \to \infty$, so that \eqref{finint} holds. 
Combining \eqref{uineq}, \eqref{finint}, and using the 
Poincar\'e inequality $\|\Delta \tilde u\|_{L^2} \ge 2\pi 
\|\nabla \tilde u\|_{L^2}$, we easily obtain
\begin{equation}\label{ubdd}
  \sup_{t \ge 0} \,e^{\mu t} \|\nabla\tilde u(t)\|_{L^2}
  \,<\, \infty~, \quad \hbox{for any } \mu < 4\pi^2~.
\end{equation}
Note, however, that the linear decay rate $\mu = 4\pi^2$ cannot be
reached by this argument, because $\int_0^\infty \|\nabla\bar 
u(t)\|_{L^3} \dd t = +\infty$ in general.

For later use, we mention that similar decay estimates 
can also be obtained for $\|\Delta\tilde u\|_{L^2}$, by 
differentiating \eqref{ns3d} and repeating the same arguments. 
We thus obtain
\begin{equation}\label{Deltabdd}
  \sup_{t \ge 1} \,e^{\mu t} \|\Delta\tilde u(t)\|_{L^2}
  \,<\, \infty~, \quad \hbox{for any } \mu < 4\pi^2~.
\end{equation}

\subsection{Evanescence of $\bar u_3$}
\label{3.2}

We next consider the third component of the two-dimensional velocity
$\bar u$, which according to \eqref{u3evol} satisfies the evolution 
equation
\begin{equation}\label{u3eq}
  \partial_t \bar u_3 + (\bar u_h \cdot \nabla)\bar u_3 + N_1
  \,=\, \Delta \bar u_3~,  
\end{equation}
where $\bar u_h = (\bar u_1,\bar u_2)^t$. The inhomogeneous term 
$N_1$ in \eqref{u3eq} is clearly negligible for large times, because
$\|N_1\|_{L^2} \le \||\tilde u||\nabla\tilde u|\|_{L^2} \le C\|\Delta
\tilde u\|_{L^2}^2$ so that $\int_0^\infty \|N_1(t)\|_{L^2}\dd t <
\infty$. By Duhamel's formula, the solution of \eqref{u3eq} can be
represented as
\begin{equation}\label{duhamel}
  \bar u_3(t) \,=\, S_{\bar u}(t,t_0)\bar u_3(t_0) -\int_{t_0}^t
  S_{\bar u}(t,s)N_1(s)\dd s~, \quad t \ge t_0 \ge 0~,
\end{equation}
where $S_{\bar u}(t,t_0)$ is the two-parameter evolution operator 
associated to the linear convection-diffusion equation $\partial_t f +
(\bar u_h \cdot \nabla)f = \Delta f$ in $\RR^2$. As is well-known
\cite{osada,carlen}, the operator $S_{\bar u}$ can be expressed by an
integral formula
\[
  (S_{\bar u}(t,t_0)f)(x) \,=\, \int_{\RR^2} \Gamma_{\bar u}
  (t,x;t_0,x_0)f(x_0)\dd x_0~,
  \quad t > t_0 \ge 0~,
\]
where the kernel $\Gamma_{\bar u}(t,x;t_0,x_0)$ has the following 
properties: 

\noindent{\bf i)} For any $\beta \in (0,1)$ there exists $C_\beta > 0$ 
such that 
\begin{equation}\label{gamm1}
  0 \,<\, \Gamma_{\bar u}(t,x;t_0,x_0) \,\le\, \frac{C_\beta}{t-t_0}
  \,\exp\Bigl(-\beta\frac{|x-x_0|^2}{4(t-t_0)}\Bigr)~,
\end{equation}
for all $t > t_0 \ge 0$ and all $x,x_0 \in \RR^2$.

\noindent{\bf ii)} For any $t > t_0 \ge 0$ and any $x,x_0 \in \RR^2$, 
one has
\begin{equation}\label{gamm2}
  \int_{\RR^2} \Gamma_{\bar u}(t,x;t_0,x_0)\dd x \,=\, 1~, \qquad 
  \int_{\RR^2} \Gamma_{\bar u}(t,x;t_0,x_0)\dd x_0 \,=\, 1~.
\end{equation}
It is very important to note that estimate \eqref{gamm1} holds 
uniformly for all $t > t_0$, with a constant $C_\beta$ which is 
independent of time. This is because $\bar\omega_3 = \partial_1
\bar u_2 - \partial_2\bar u_1$ is uniformly bounded in $L^1(\RR^2)$, 
see \cite{osada}. It follows in particular from \eqref{gamm1}, 
\eqref{gamm2} that $\|S_{\bar u}(t,t_0)f\|_{L^2} \le \|f\|_{L^2}$ 
for all $t \ge t_0$, and that $S_{\bar u}(t,t_0)$ satisfies similar
$L^p$--$L^q$ estimates as the heat semigroup $e^{(t-t_0)\Delta}$. 

We claim that the solution $\bar u_3(t)$ of \eqref{u3eq} converges 
to zero in $L^2(\RR^2)$ as $t \to \infty$. To prove that, fix any
$\epsilon > 0$, and take $t_0 > 0$ sufficiently large so that 
$\int_{t_0}^\infty \|N_1(s)\|_{L^2}\dd s \le \epsilon$. Then
\[
  \Bigl\|\int_{t_0}^t S_{\bar u}(t,s)N_1(s)\dd s\Bigr\|_{L^2}
  \,\le\, \int_{t_0}^\infty \|N_1(s)\|_{L^2}\dd s \,\le\, \epsilon~,
  \quad \hbox{for all } t \ge t_0~,
\]
hence in the right-hand side of \eqref{duhamel} it is sufficient to
bound the first term $v(t) = S_{\bar u}(t,t_0)\bar u_3(t_0)$.  Since
$\bar u_3(t_0) \in L^2(\RR^2)$, we can decompose $\bar u_3(t_0) = v_1 +
v_2$ with $v_1 \in L^1(\RR^2) \cap L^2(\RR^2)$ and $\|v_2\|_{L^2}
\le \epsilon$. Then $v(t) = v_1(t) + v_2(t)$ with
\[
  \|v_1(t)\|_{L^2} \,=\, \|S_{\bar u}(t,t_0)v_1\|_{L^2}
  \,\le\, \frac{C}{(t-t_0)^{1/2}}\,\|v_1\|_{L^2} \,\xrightarrow[t 
  \to \infty]{}\, 0~,
\]
and $\|v_2(t)\|_{L^2} = \|S_{\bar u}(t,t_0)v_2\|_{L^2} \le
\|v_2\|_{L^2} \le \epsilon$. Thus, if $t > t_0$ is sufficiently large, 
we have
\[
  \|\bar u_3(t)\|_{L^2} \,\le\, \|v_1(t)\|_{L^2} + \|v_2(t)\|_{L^2}
  + \int_{t_0}^t \|N_1(s)\|_{L^2}\dd s \,\le\, 3\epsilon~,
\]
which proves the claim. 

On the other hand, we know from \eqref{last4} that $\int_0^\infty 
\|\nabla \bar u_3(t)\|_{L^2}^2\dd t < \infty$, hence there exists
a sequence $t_n \to \infty$ such that $\|\nabla \bar u_3(t_n)
\|_{L^2}^2 \to 0$ as $n \to \infty$. In view of \eqref{sys2}, 
we have for each $n$:
\[
  \sup_{t \ge t_n}\|\nabla \bar u_3(t)\|_{L^2}^2 \,\le\, 
  \|\nabla \bar u_3(t_n)\|_{L^2}^2 + C\int_{t_n}^\infty
  (\|\nabla \bar u_3(s)\|_{L^2}^2 \|\bar\omega_3(t)\|_{L^2}^2
  + \|\nabla\tilde u(s)\|_{L^2}^3 \|\Delta \tilde u(s)\|_{L^2}
  )\dd s~, 
\]
and the right-hand side converges to zero as $n \to \infty$. 
This shows that $\|\nabla \bar u_3(t)\|_{L^2} \to 0$ as 
$t \to \infty$, and we have therefore proved that $\bar u_3(t)$
converges to zero in $H^1(\RR^2)$ as $t \to \infty$. 

\subsection{Diffusive estimates for $\bar\omega_3$}
\label{3.3}

We now turn our attention to the third component of the
two-dimensional vorticity $\bar \omega$, which evolves 
according to \eqref{omega3evol}:
\begin{equation}
\label{omega3eq}
  \partial_t \bar \omega_3 + (\bar u_h \cdot \nabla)\bar \omega_3 
  + N_2 \,=\, \Delta \bar \omega_3~.  
\end{equation}
By \eqref{last4}, there exists $C_4 > 0$ such that $\|\bar\omega_3(t)
\|_{L^1} + \|\bar\omega_3(t)\|_{L^2} \le C_4$ for all $t \ge 0$.  To
obtain sharper estimates, including decay rates in time, we use a
standard method that goes back to Nash, see \cite{fabes}. By the
Gagliardo-Nirenberg inequality, there exists $C > 0$ such that
$\|\bar\omega_3\|_{L^2}^2 \le C \|\bar\omega_3\|_{L^1}\|\nabla\bar
\omega_3\|_{L^2}$, hence $\|\bar\omega_3\|_{L^2}^2 \le C C_4
\|\nabla\bar\omega_3\|_{L^2}$. Inserting this bound into \eqref{sys3}, 
we obtain
\begin{equation}\label{omeganash}
  \frac{\dD}{\dD t}\|\bar \omega_3(t)\|_{L^2}^2 \,\le\, 
  - C_5 \|\bar \omega_3(t)\|_{L^2}^4 + 8 \||\tilde u(t)| |\nabla
  \tilde u(t)|\|_{L^2}~,
\end{equation}
where $C_5 = (CC_4)^{-2}$. Since $\||\tilde u(t)||\nabla\tilde u(t)|
\|_{L^2}$ decays exponentially to zero as $t \to \infty$, it 
follows from \eqref{omeganash} that
\begin{equation}\label{L2tbound}
  \sup_{t \ge 0}(1+t) \|\bar \omega_3(t)\|_{L^2}^2 \,=\, C_6 \,<\,
  \infty~.
\end{equation}

A similar argument can be used to estimate $\|\nabla\bar \omega_3
\|_{L^2}$. From \eqref{omega3eq} we have
\[
  \frac12 \frac{\dD}{\dD t}\|\nabla\bar \omega_3\|_{L^2}^2 
  \,=\, -\int_{\RR^2} |\Delta\bar\omega_3|^2\dd x + 
  \int_{\RR^2} (\Delta\bar\omega_3)(\bar u_h \cdot \nabla)\bar 
  \omega_3 \dd x + \int_{\RR^2} (\Delta\bar\omega_3)N_2\dd x~.
\]
Integrating by parts and using the fact that $\|\nabla\bar 
u_h\|_{L^2} = \|\bar\omega_3\|_{L^2}$, we find
\begin{align*}
  \Bigl|\int_{\RR^2} (\Delta\bar\omega_3)(\bar u_h \cdot \nabla)\bar 
  \omega_3 \dd x\Bigr| \,&\le\, \| |\nabla\bar\omega_3|
  \,|\nabla\bar u_h|\,|\nabla\bar\omega_3|\|_{L^1} \,\le\, 
  \|\nabla\bar\omega_3\|_{L^4}^2 \|\bar\omega_3\|_{L^2} \\
  \,&\le\, C \|\Delta\bar\omega_3\|_{L^2} \|\nabla\bar\omega_3\|_{L^2}
  \|\bar\omega_3\|_{L^2} \,\le\, C \|\Delta\bar\omega_3\|_{L^2}^{3/2} 
  \|\bar\omega_3\|_{L^2}^{3/2}~,
\end{align*}
hence
\[
  \frac{\dD}{\dD t}\|\nabla\bar \omega_3(t)\|_{L^2}^2 
  \,\le\, -\|\Delta\bar\omega_3(t)\|_{L^2}^2 + C(\|\bar\omega_3(t)
  \|_{L^2}^6 + \|N_2(t)\|_{L^2}^2)~.
\]
As $\|\nabla\bar \omega_3\|_{L^2}^2 \le \|\bar \omega_3\|_{L^2}
\|\Delta\bar \omega_3\|_{L^2} \le C_6^{1/2}(1+t)^{-1/2}
\|\Delta\bar \omega_3\|_{L^2}$, we conclude that
\begin{equation}\label{nabla3ineq}
  \frac{\dD}{\dD t}\|\nabla\bar \omega_3(t)\|_{L^2}^2 
  \,\le\, -C_6^{-1}(1+t)\|\nabla\bar\omega_3(t)\|_{L^2}^4 + 
  C(\|\bar\omega_3(t)\|_{L^2}^6 + \|N_2(t)\|_{L^2}^2)~.
\end{equation}
Now, since $\|\bar\omega_3(t)\|_{L^2}^6 \le C_6^3(1+t)^{-3}$, and since
$\|N_2(t)\|_{L^2}^2$ decays exponentially to zero as $t \to \infty$, 
the differential inequality \eqref{nabla3ineq} implies that 
$\|\nabla\bar \omega_3(t)\|_{L^2}^2$ decreases at least like $t^{-2}$ 
as $t \to \infty$. Taking into account the fact that $\bar
\omega_3(0) \in L^2(\RR^2)$, we arrive at
\begin{equation}  \label{H1tbound}
  \sup_{t \ge 0}t(1+t) \|\nabla\bar\omega_3(t)\|_{L^2}^2 \,=\, 
  C_7 \,<\,\infty~.
\end{equation}

\subsection{Compactness of the rescaled solution}
\label{3.4}

To show that the solution $\bar\omega_3(t,x)$ of \eqref{omega3eq}
converges to Oseen's vortex as $t \to \infty$, it is convenient to
introduce self-similar variables. Following \cite{thgcew,thgcew2}, 
we define
\begin{align}\label{scaling}
  &\bar\omega_3(t,x) \,=\,\frac{1}{1+t} \,w\left(\log(1+t)\,,\,
  \frac{x}{\sqrt{1+t}}\right)~, \\ \nonumber
  &\bar u_h(t,x) \,=\,\frac{1}{\sqrt{1+t}} \,v\left(\log(1+t)\,,\,
  \frac{x}{\sqrt{1+t}}\right)~.
\end{align}
We also denote
\[
  \xi \,=\, \frac{x}{\sqrt{1+t}}~, \qquad \tau \,=\, \log(1+t)~.
\]
Then the rescaled vorticity $w(\tau,\xi)$ satisfies the equation
\begin{equation}\label{weq}
  \partial_\tau w + (v\cdot \nabla_\xi) w + \tilde N_2  
  \,=\, \Delta_\xi w + \frac12 (\xi\cdot\nabla_\xi) w + w~,
\end{equation}
where $\tilde N_2(\tau,\xi) = e^{2\tau}N_2(e^\tau{-}1,\xi\,e^{\tau/2})$, 
and $v(\tau,\xi)$ coincides with the two-dimensional velocity field 
obtained from $w(\tau,\xi)$ via the Biot-Savart law \eqref{BSbar2}. 
It is clear that
\[
  \int_0^\infty \|\tilde N_2(\tau)\|_{L^1}\dd\tau \,=\, 
  \int_0^\infty e^\tau \|N_2(e^\tau-1)\|_{L^1}\dd\tau \,=\, 
  \int_0^\infty \|N_2(t)\|_{L^1}\dd t \,<\, \infty~,
\]
hence the term $\tilde N_2(\tau,\xi)$ in \eqref{weq} will be negligible 
for large times. The solution of \eqref{weq} can be represented 
as 
\begin{equation}\label{duhamel2}
  w(\tau) \,=\, \tilde S_v(\tau,\tau_0)w(\tau_0) -\int_{\tau_0}^\tau
  \tilde S_v(\tau,s)\tilde N_2(s)\dd s~, \quad \tau \ge \tau_0 \ge 0~,
\end{equation}
where in analogy with \eqref{duhamel} we denote by $\tilde 
S_v(\tau,\tau_0)$ the two-parameter evolution operator associated 
to the linear equation $\partial_\tau w + (v\cdot\nabla)w = 
\Delta w + \frac12 (\xi\cdot\nabla) w + w$ (note that $\tilde S_v$ 
depends on the velocity field $v(\tau,\xi)$, which is considered 
here as given). Using the same notations as in Section~\ref{3.2}, 
we find that
\begin{equation}\label{Gammau}
  (\tilde S_v(\tau,\tau_0)f)(\xi) \,=\, \int_{\RR^2} e^\tau
  \Gamma_{\bar u}(e^\tau-1,\xi\,e^{\tau/2};e^{\tau_0}-1,\xi_0\,
  e^{\tau_0/2})f(\xi_0)\dd \xi_0~.
\end{equation}

The aim of this paragraph is to prove the following basic result:

\begin{lm}\label{compact}
The solution $\{w(\tau)\}_{\tau\ge 0}$ of \eqref{weq} is relatively 
compact in $L^1(\RR^2)$. 
\end{lm}

\noindent\textbf{Proof.} By construction $w \in C^0([0,\infty),
L^1(\RR^2))$ and $\|w(\tau)\|_{L^1} \le C_4$ for all $\tau \ge 0$. 
To prove compactness, we use the Riesz criterion \cite{reedsimon} 
and proceed in two steps:

\noindent{\bf i)} We first show that
\begin{equation}\label{loc1}
  \sup_{\tau\ge 0}\int_{|\xi|\ge R}|w(\tau,\xi)|\dd \xi
  \,\xrightarrow[R \to \infty]{}\, 0~. 
\end{equation}
Indeed, fix $\epsilon > 0$ and take $\tau_0 \ge 0$ large enough 
so that $\int_{\tau_0}^\infty \|\tilde N_2(\tau)\|_{L^1}\dd \tau
\le \epsilon/2$. Then choose $R_1 \ge 0$ large enough so that
\[
  \sup_{\tau\in[0,\tau_0]}\int_{|\xi|\ge R_1}|w(\tau,\xi)|\dd \xi
  \,\le\, \epsilon~.
\]
This is clearly possible, because the finite-time trajectory 
$\{w(\tau)\,|\, 0 \le \tau \le \tau_0\}$ is compact in 
$L^1(\RR^2)$. For $\tau \ge \tau_0$ the solution of \eqref{weq}
can be represented as in \eqref{duhamel2}, where the second 
term in the right-hand side satisfies
\[
  \Bigl\|\int_{\tau_0}^\tau \tilde S_v(\tau,s)\tilde N_2(s)\dd s 
  \Bigr\|_{L^1} \,\le\, \int_{\tau_0}^\tau \|\tilde N_2(s)\|_{L^1}
  \dd s \,\le\, \epsilon/2~. 
\]
As for the first term $w_1(\tau) = \tilde S_v(\tau,\tau_0)w(\tau_0)$, 
it can be estimated by a direct calculation, using the representation
formula \eqref{Gammau} and the bounds \eqref{gamm1} on the kernel
$\Gamma_{\bar u}$. Proceeding exactly as in the proof of 
\cite[Lemma~2.5]{thgcew2}, one finds $R_2 \ge 0$ such that 
\[
  \sup_{\tau \ge \tau_0}\int_{|\xi|\ge R_2}|w_1(\tau,\xi)|\dd \xi
  \,\le\, \frac{\epsilon}{2}~.
\]
If we now choose $R = \max(R_1,R_2)$, we see that $\int_{|\xi|\ge R}|
w(\tau,\xi)|\dd \xi \le \epsilon$ for all $\tau \ge 0$, which 
proves \eqref{loc1}. 

\noindent{\bf ii)} Our second task is to verify that
\begin{equation}\label{loc2}
  \sup_{\tau \ge 0} \,\sup_{|\eta| \le \delta}\int_{\RR^2}
  |w(\tau,\xi-\eta) - w(\tau,\xi)|\dd \xi 
  \,\xrightarrow[\delta \to 0]{}\,0~.
\end{equation}
By compactness of the finite-time trajectory, it is sufficient
to check \eqref{loc2} for $\tau \ge 1$. Using the definitions
\eqref{scaling} and the bound \eqref{H1tbound} established in 
Section~\ref{3.3}, we find
\[
  \sup_{\tau \ge 1} \|\nabla w(\tau)\|_{L^2} \,=\, C_8 \,<\, \infty~.
\]
Fix $\epsilon > 0$. By the first step, there exists $R \ge 1$ 
such that
\[
  \sup_{\tau \ge 1}\int_{|\xi|\ge R-1}|w(\tau,\xi)|\dd \xi
  \,\le\, \frac{\epsilon}{3}~.
\]
Take $\delta \in (0,1]$ such that $C_8 \delta \pi^{1/2}(R+1) \le
\epsilon/3$. If $\eta \in \RR^2$ satisfies $|\eta| \le \delta$, 
we have
\[
  \int_{|\xi| \ge R} |w(\tau,\xi-\eta) - w(\tau,\xi)|\dd \xi \,\le\, 
  2 \int_{|\xi| \ge R-1} |w(\tau,\xi)|\dd \xi \,\le\, 
  \frac{2\epsilon}{3}~.
\]
On the other hand, by Fubini's theorem and H\"older's inequality, 
\begin{align*}
  \int_{|\xi| \le R} |w(\tau,\xi-\eta) - w(\tau,\xi)|\dd \xi \,&\le\, 
  \int_{|\xi| \le R} \int_0^1 |\eta\cdot\nabla w(\tau,\xi-r\eta)|\dd r
  \dd \xi \\
  \,&\le\, |\eta| \int_{|\xi| \le R+1} |\nabla w(\tau,\xi)|\dd\xi
  \,\le\, C_8 |\eta| \pi^{1/2}(R+1) \,\le\, \frac{\epsilon}{3}~,
\end{align*}
hence $\int_{\RR^2} |w(\tau,\xi-\eta) - w(\tau,\xi)|\dd \xi \le
\epsilon$ for all $\tau \ge 1$ whenever $|\eta| \le \delta$. 
This proves \eqref{loc2}. By the Riesz criterion, \eqref{loc1}
and \eqref{loc2} together imply that the trajectory 
$\{w(\tau)\}_{\tau\ge0}$ is relatively compact in $L^1(\RR^2)$. 
\rule{2mm}{2mm}

\subsection{Determination of the $\omega$-limit set}
\label{3.5}

We know from Lemma~\ref{compact} that the solution 
$\{w(\tau)\}_{\tau\ge 0}$ of \eqref{weq} lies in a compact subset
of $L^1(\RR^2)$. Let $\Omega_\infty$ be the $\omega$-limit set 
of this solution, namely
\[
  \Omega_\infty \,=\, \Bigl\{w_\infty \in L^1(\RR^2)\,\Big|\,
  \exists \tau_n \to \infty \hbox{ such that } w(\tau_n) 
  \xrightarrow[n \to \infty]{L^1} w_\infty\Bigr\}~.
\]
Since $\int_{\RR^2} w(\tau,\xi)\dd\xi = \int_{\RR^2} \bar
\omega_3(e^\tau{-}1,x)\dd x = \alpha$ for all $\tau \ge 0$, where
$\alpha$ is given by \eqref{alphadef}, it is clear that 
\begin{equation}\label{intcond}
  \int_{\RR^2} w_\infty(\xi)\dd\xi \,=\, \alpha~, \quad 
  \hbox{for all } w_\infty \in \Omega_\infty~.  
\end{equation}
Our goal is to show that $\Omega_\infty = \{\alpha g\}$,
where $g(\xi) = (4\pi)^{-1} e^{-|\xi|^2/4}$. This will imply that
$\|w(\tau)-\alpha g\|_{L^1} \to 0$ as $\tau \to \infty$, which is
equivalent to \eqref{omegaconv}.

Let $\hat\Phi(\tau)$ denote the semiflow defined by the limiting 
equation
\begin{equation}\label{weq2}
  \partial_\tau \hat w + \hat v\cdot \nabla_\xi \hat w \,=\, \Delta_\xi 
  \hat w + \frac12 \xi\cdot\nabla_\xi \hat w + \hat w~,
\end{equation}
where $\hat v$ is the velocity field obtained from $\hat w$ via the 
Biot-Savart law \eqref{BSbar2}. Note that \eqref{weq2} is just 
the ordinary two-dimensional vorticity equation expressed in 
self-similar variables. We shall prove that the $\omega$-limit
set of the solution $w(\tau)$ of \eqref{weq} is totally 
invariant under the evolution defined by \eqref{weq2}: 

\begin{lm}\label{invariant}
The $\omega$-limit set $\Omega_\infty$ satisfies $\hat\Phi(\tau)
\Omega_\infty = \Omega_\infty$ for all $\tau \ge 0$. 
\end{lm}

Using \cite[Proposition~3.5]{thgcew2}, we deduce that $\Omega_\infty
\subset \{\alpha' g\,|\,\alpha' \in \RR\}$, hence $\Omega_\infty
 = \{\alpha g\}$ in view of \eqref{intcond}. This is the desired
result, which completes the proof of Theorem~\ref{main}. 

\medskip
\noindent\textbf{Proof of Lemma~\ref{invariant}.} Let 
$\{S(\tau)\}_{\tau \ge 0}$ denote the $C_0$-semigroup generated by the
Fokker-Planck operator $\Delta + \frac12 \xi\cdot\nabla + 1$, see
\cite{thgcew}. If $w \in L^1(\RR^2)$, then for any $p \in [1,\infty)$
we have the following estimates:
\begin{equation}\label{l1lp}
  \|S(\tau)w\|_{L^p} \,\le\, \frac{\|w\|_{L^1}}{4\pi 
  a(\tau)^{1-\frac{1}{p}}}~, \qquad
  \|\nabla S(\tau)w\|_{L^p} \,\le\, \frac{C\|w\|_{L^1}}{ 
  a(\tau)^{\frac32-\frac{1}{p}}}~, \quad \tau > 0~,
\end{equation}
where $a(\tau) = 1-e^{-\tau}$. Moreover $\|S(\tau)w\|_{L^p} 
\le e^{\tau(1-\frac1p)}\|w\|_{L^p}$ for all $\tau \ge 0$ if
$w \in L^p(\RR^2)$. 

Let $w_\infty \in \Omega_\infty$, and take a sequence $\tau_n \to \infty$ 
such that $\|w(\tau_n)-w_\infty\|_{L^1} \to 0$ as $n \to \infty$. 
Since the trajectory $\{w(\tau)\}_{\tau\ge 0}$ is bounded in $L^2(\RR^2)$ 
by \eqref{L2tbound}, \eqref{scaling}, we have $w_\infty \in 
L^2(\RR^2)$ and (up to extracting a subsequence) we can assume
that $\|w(\tau_n)-w_\infty\|_{L^p} \to 0$ as $n \to \infty$
for any $p \in [1,2)$. For each $n \in \NN$, let $w_n(\tau) = 
w(\tau+\tau_n)$ and $v_n(\tau) = v(\tau+\tau_n)$. Then $w_n(\tau)$ 
satisfies the integral equation
\begin{equation}\label{wneq}
  w_n(\tau) \,=\, S(\tau)w(\tau_n) -\int_0^\tau S(\tau-s)
  \Bigl(v_n(s)\cdot\nabla w_n(s) + \tilde N_2(\tau_n+s)\Bigr)
  \dd s~.
\end{equation}
On the other hand, if we denote $\hat w(\tau) = \hat\Phi(\tau)w_\infty$, 
we have
\begin{equation}\label{hatweq}
  \hat w(\tau) \,=\, S(\tau)w_\infty -\int_0^\tau S(\tau-s)
  \hat v(s)\cdot\nabla \hat w(s) \dd s~.
\end{equation}
Subtracting \eqref{hatweq} from \eqref{wneq} and using the bounds
\eqref{l1lp} on the semigroup $S(\tau)$, we obtain for any 
$p \in [1,2)$:
\begin{align}\nonumber
  \|w_n(\tau)&-\hat w(\tau)\|_{L^p} \,\le\, 
  e^{\tau(1-\frac1p)}\|w(\tau_n)-w_\infty\|_{L^p} + \int_0^\tau \frac{C}{
  a(\tau-s)^{1-\frac1p}}\,\|\tilde N_2(\tau_n+s)\|_{L^1}\dd s \\
  \,&+\,\int_0^\tau \frac{C\,e^{-\frac12(\tau-s)}}{a(\tau-s)^{
  \frac32-\frac1p}}(\|w_n(s)\|_{L^{4/3}} + \|\hat w(s)\|_{L^{4/3}}) 
  \|w_n(s)-\hat w(s)\|_{L^{4/3}}\dd s~.\label{diffp}
\end{align}
Here we have used the fact that $S(\tau)v\cdot\nabla w =
S(\tau)\nabla\cdot(vw) = e^{-\tau/2}\nabla\cdot S(\tau)(vw)$, and the
bound $\|vw\|_{L^1} \le \|v\|_{L^4}\|w\|_{L^{4/3}} \le
C\|w\|_{L^{4/3}}^2$ which holds in view of Proposition~\ref{BSbarprop}.
We first choose $p = 4/3$ and consider equation \eqref{diffp} for
$\tau$ in some compact interval $[0,T]$. The first line in the
right-hand side converges uniformly to zero as $n \to \infty$, and in
the second line we know that $\|w_n(s)\|_{L^{4/3}} + \|\hat
w(s)\|_{L^{4/3}}$ is uniformly bounded for all $n \in \NN$ and all
$\tau \in [0,T]$. Thus it follows from Gronwall's lemma \cite{henry}
that
\begin{equation}\label{L43conv}
  \sup_{\tau\in[0,T]} \|w_n(\tau) - \hat w(\tau)\|_{L^{4/3}}
  \,\xrightarrow[n\to\infty]{}\,0~.
\end{equation}
Setting now $p=1$ in \eqref{diffp} and using \eqref{L43conv}, 
we conclude that $\|w_n(\tau) - \hat w(\tau)\|_{L^1} \to 0$ 
as $n \to \infty$, for all $\tau \in [0,T]$. In other words 
$w(\tau+\tau_n)$ converges to $\hat\Phi(\tau)w_\infty$ as $n \to \infty$, 
which means that $\hat\Phi(\tau)w_\infty \in \Omega_\infty$ for all 
$\tau \in [0,T]$. Since $T > 0$ was arbitrary, we have shown that 
$\hat\Phi(\tau)\Omega_\infty \subset \Omega_\infty$ for all $\tau \ge 0$. 

To prove the converse inclusion, we fix $\tau \ge 0$ and take again
$w_\infty \in \Omega_\infty$. If $w(\tau_n) \to w_\infty$ in
$L^1(\RR^2)$ as $n \to \infty$, then after extracting a subsequence we
can assume that $w(\tau_n-\tau)$ converges as $n \to \infty$ to
some $w_0 \in \Omega_\infty$. Using exactly the same arguments
as before, we can prove that $w_\infty = \hat\Phi(\tau)w_0$. 
This shows that $\Omega_\infty \subset \hat\Phi(\tau)\Omega_\infty$, 
for any $\tau \ge 0$. \rule{2mm}{2mm}

\appendix  
\section{Appendix : The Biot-Savart Law in $\RR^2 \times 
\TT^1$}
\label{BSlaw}

In this appendix we give explicit formulas for the Biot-Savart law 
in the domain $\DD = \RR^2 \times \TT^1$, and we collect a few
estimates for the velocity field $u$ in terms of the vorticity 
$\omega$ which are used throughout the paper. All these results
are well-known (see \cite{roussier}) and are reproduced here 
for the reader's convenience. 

Let $u : \DD \to \RR^3$ be a divergence-free velocity field, 
and denote by $\omega = \curl u$ the associated vorticity 
field. As is explained in the introduction, it is convenient
to decompose
\[
  u(x,z) \,=\, \bar{u}(x) + \tilde{u}(x,z)~, \quad 
  \omega(x,z) \,=\, \bar{\omega}(x) + \tilde{\omega}(x,z)~,
  \quad x \in \RR^2~, \quad z \in \TT^1~,
\]
where $\bar u = Qu$, $\bar\omega = Q\omega$, and $Q$ is the vertical 
average operator defined by \eqref{Qdef}. Then it is straightforward 
to verify that $\bar{\omega}=\curl{\bar{u}}$ and $\tilde{\omega}
=\curl{\tilde{u}}$. Moreover, the four vector fields $\bar u$, 
$\tilde u$, $\bar \omega$, $\tilde \omega$ are all divergence-free.
Thus we can consider separately the Biot-Savart law for the 
two-dimensional part $(\bar u,\bar \omega)$ and for the 
three-dimensional fluctuation $(\tilde u,\tilde \omega)$. 

\subsection{The Biot-Savart law for $(\bar{u},\bar{\omega})$.}  
\label{appendix A.1}  
  
Since the vector fields $\bar u$, $\bar \omega$ do not 
depend on the vertical variable $z$, the relations 
$\div \bar u = 0$ and $\curl \bar u = \bar\omega$ can
be written in the following equivalent form: 
\begin{equation}\label{absys}
  (a)\,\left\{\!\begin{array}{l}  
  \bar{\omega}_1= \partial_2 \bar{u}_3~, \\  
  \bar{\omega}_2= -\partial_1 \bar{u}_3~,
  \end{array}\right. \qquad 
  (b)\,\left\{\!\begin{array}{l}  
  \partial_1 \bar{u}_2 - \partial_2 \bar{u}_1 = \bar{\omega}_3~,\\  
  \partial_1 \bar{u}_1 + \partial_2 \bar{u}_2 = 0~.  
  \end{array}\right.
\end{equation}
To solve the first system $(a)$, we observe that $\Delta\bar u_3
= \partial_2\bar \omega_1 - \partial_1\bar \omega_2$ and we use
the fundamental solution of the Laplacian operator in $\RR^2$. 
After integrating by parts, we obtain
\begin{equation}\label{BSbar1} 
  \bar{u}_3 (x) \,=\, - \frac{1}{2\pi} \int_{\RR^2}
  \frac{(x-y)}{|x-y|^2}\wedge \begin{pmatrix} \bar{\omega}_1 
  \\ \bar{\omega}_2\end{pmatrix}(y)\dd y~, \quad x \in \RR^2~.  
\end{equation}  
On the other hand, the solution of system $(b)$ is just the ordinary 
Biot-Savart law in $\RR^2$:
\begin{equation}\label{BSbar2}  
  \begin{pmatrix}\bar{u}_1 \\ \bar{u}_2   
  \end{pmatrix}(x) \,=\, \frac{1}{2\pi} \int_{\RR^2} 
  \frac{(x-y)^{\perp}}{|x-y|^2}\,\bar{\omega}_3(y)\dd y~, \quad 
  x \in \RR^2~.
\end{equation}
Here, if $x=(x_1,x_2) \in \RR^2$, we denote $x^{\perp}=(-x_2,x_1)$. 
In particular, we see from \eqref{BSbar2} that the horizontal part $\bar
u_h = (\bar u_1,\bar u_2)$ of the velocity field $\bar u$ can be
reconstructed from the third component $\bar\omega_3$ of the vorticity
$\bar \omega$, an observation that is used many times in the
previous sections. 

In both formulas \eqref{BSbar1} and \eqref{BSbar2}, the velocity
field is expressed in terms of the vorticity through a convolution
with a singular integral kernel, which is homogeneous of degree
$-1$. Thus we can apply the classical Hardy-Littlewood-Sobolev 
inequality \cite{liebloss} to both cases, and obtain the 
following result:
  
\begin{prop}\label{BSbarprop} 
Let $\bar{u}$ be the velocity field obtained from $\bar{\omega}$ via
the Biot-Savart law \eqref{BSbar1}, \eqref{BSbar2}. Assume that
$1< p < 2 < q< \infty$ and $\frac{1}{q} = \frac{1}{p} - \frac{1}{2}$. 
If $\bar{\omega} \in L^p(\RR^2)^3$, then $\bar{u}\in L^q(\RR^2)^3$, and 
there exists a constant $C > 0$ (depending only on $p$) such that
\[
  \|\bar{u}\|_{L^q(\RR^2)} \,\le\, C \|\bar{\omega}\|_{L^p(\RR^2)}~.
\]
\end{prop}  

Moreover, using Calder\'on-Zygmund's theory, one can show that
$\|\nabla \bar u\|_{L^p} \le C\|\bar\omega\|_{L^p}$ for $1 < p < \infty$. 
In the particular case $p = 2$, we even have $\|\nabla \bar u\|_{L^2} 
= \|\bar\omega\|_{L^2}$. 

\subsection{The Biot-Savart law for $(\tilde{u},\tilde{\omega})$}  
\label{appendix A.2}  

The relation between $\tilde u$ and $\tilde \omega$ is most
conveniently expressed in Fourier variables. Using the same notations 
as in \eqref{fourier1}, we can write
\begin{equation}\label{uofourier}
  \tilde u(x,z) \,=\, \int_{\RR^2} \sum_{n \in \ZZ^*} \tilde u_n(k) 
  \,e^{i(k\cdot x +2\pi nz)}\,\frac{\dD k}{2\pi}~, \quad
  \tilde \omega(x,z) \,=\, \int_{\RR^2} \sum_{n \in \ZZ^*} \tilde 
  \omega_n(k)\,e^{i(k\cdot x +2\pi nz)}\,\frac{\dD k}{2\pi}~.
\end{equation}
Observe that the sums here are taken over $n \in \ZZ^* \equiv 
\ZZ\setminus\{0\}$, because $\tilde u$ and $\tilde \omega$ have zero 
average with respect to the vertical variable. Since $\div\tilde u 
= 0$ and $\curl \tilde u = \tilde \omega$, we have $-\Delta\tilde u 
= \curl \tilde\omega$, hence
\begin{equation}\label{BStilde}  
  \tilde{u}_n(k) \,=\, \frac{1}{|k|^2+4\pi^2 n^2} 
  \begin{pmatrix} 0 & -2\pi in & i k_2 \\  
  2\pi in & 0 & -i k_1 \\  -ik_2 & ik_1 & 0  
  \end{pmatrix}\tilde{\omega}_n(k)~, \quad n \in \ZZ^*~, 
  \quad k \in \RR^2~.  
\end{equation}  
Since $n \neq 0$ in \eqref{BStilde}, it follows that 
$\|\tilde{u}\|_{H^{s+1}} \le C \|\tilde{\omega}\|_{H^s}$ 
for any $s \ge 0$, see \eqref{Hsnorm}. In particular, taking 
$s = 0$ and using the Sobolev embedding $H^1(\DD) \hookrightarrow
L^q(\DD)$ for $q \in [2,6]$, we obtain: 

\begin{prop}  
\label{BStildeprop}
Let $\tilde{u}$ be the velocity field obtained from $\tilde{\omega}$
via the Biot-Savart law \eqref{BStilde}. If $\tilde{\omega} \in
L^2(\DD)$, then $\tilde{u} \in L^q(\DD)$ for any $q \in [2,6]$, 
and there exists $C>0$ (depending only on $q$) such that 
\[
  \|\tilde{u}\|_{L^q(\DD)} \,\le\, C \|\tilde{\omega}\|_{L^2(\DD)}~.
\]
\end{prop}  

\subsection{The Leray projector}  
\label{appendix A.3}  

In the Fourier variables defined by \eqref{fourier1},
\eqref{fourier2}, the Leray projector $\PP$ has the following 
simple expression
\begin{equation}\label{PPdef}
  (\PP f)_n(k) \,=\, f_n(k) + \frac{\xi \cdot f_n(k)}{|k|^2 
  +4\pi^2 n^2}\,\xi~, \quad \hbox{where} \quad \xi \,=\, 
  \begin{pmatrix} ik\\2 \pi in \end{pmatrix} \in \RR^3~.
\end{equation}
Clearly $\PP$ commutes with the vertical average operator 
$Q$, which satisfies $(Qf)_n(k) \,=\, f_n(k)\delta_{n,0}$. If 
$\bar f = Qf$, we see from \eqref{PPdef}
that $e_3 \cdot (\PP \bar f) = e_3 \cdot \bar f$. In other 
words, the Leray projector $\PP$ acts trivially on the third 
component of $z$-independent vector fields. 

\section{Appendix: Dispersive estimates}
\label{dispersive}  

This final section is devoted to the proof of
Proposition~\ref{strichartz}.  The arguments here follow closely the
analysis of \cite[Chap.~5]{cdgg}, and were already published 
in \cite{roussier2} in a slightly different form. 

\medskip
\noindent{\bf Proof of proposition \ref{strichartz}:}
If $\tilde u(t,x,z)$ is a divergence-free solution of the linear 
Rossby equation \eqref{freerot}, we first observe that the 
Fourier transform $\tilde u_n(t,k)$, which is defined as in 
\eqref{uofourier}, satisfies
\[
  \partial_t \tilde u_n(t,k) + M_n^\Omega(k)\tilde u_n(t,k) \,=\, 0~,
  \quad k \in \RR^2~, \quad n \in \ZZ^*~,
\]
where $M_n^\Omega(k)$ is the $3\times3$ matrix defined by
\begin{equation}\label{Ankdef}
  M_n^{\Omega}(k) \,=\, (|k|^2 +4\pi^2n^2)\one + \frac{2i\pi n 
  \Omega}{|k|^2 +4\pi^2n^2} \begin{pmatrix} 0 & 
  -2\pi in & ik_2 \\2\pi in & 0 & -ik_1 \\ -ik_2 & ik_1 & 0 
  \end{pmatrix}~.
\end{equation}
Indeed, the first term in \eqref{Ankdef} corresponds to 
$-(\Delta\tilde u)_n(k) = (|k|^2 +4\pi^2n^2)\tilde u_n(k)$.
On the other hand, if $\tilde\omega = \curl\tilde u$, we have
from \eqref{BStilde}
\[
  \tilde \omega_n(k) \,=\, \xi \wedge \tilde u_n(k)~, \quad
  \tilde u_n(k) \,=\, \frac{\xi \wedge \tilde \omega_n(k)}{|\xi|^2}~, 
  \quad \hbox{where}\quad \xi \,=\, \begin{pmatrix} ik\\2 \pi in 
  \end{pmatrix}~.
\]
It follows that
\[
  e_3 \wedge \tilde u_n(k) \,=\, \frac{1}{|\xi|^2}\,
  e_3 \wedge (\xi \wedge \tilde \omega_n(k)) \,=\,
  \frac{1}{|\xi|^2}\,\Bigl((e_3\cdot\tilde\omega_n(k))\xi
  - (e_3\cdot\xi)\tilde\omega_n(k)\Bigr)~.
\]
The last member is the sum of two terms, one of which is 
proportional to $\xi$ (gradient term) and the other orthogonal 
to $\xi$ (divergence-free term). Thus
\[
 -\PP(e_3 \wedge \tilde u_n(k)) \,=\, \frac{1}{|\xi|^2}\,
 (e_3\cdot\xi)\tilde\omega_n(k) \,=\, \frac{2\pi in}{|\xi|^2}
 \,\xi\wedge\tilde u_n(k)~,
\]
which gives the second term in \eqref{Ankdef}. 

As is easily verified, the eigenvalues of $M_n^\Omega(k)$ are 
$|\xi|^2$ and $|\xi|^2 \pm i\Omega \eta$, 
where
\begin{equation}\label{xietadef}
  |\xi| \,=\, |\xi(k,n)| \,=\, \sqrt{|k|^2+4\pi^2 n^2}~, \quad\hbox{and}
  \quad \eta  \,=\, \eta(k,n)  \,=\, \frac{2\pi n}{\sqrt{|k|^2+4\pi^2 n^2}}~.
\end{equation}
Moreover, the eigenvector corresponding to $|\xi|^2$ is 
proportional to $\xi$, whereas the normalized eigenvectors $w_n^\pm(k)$ 
corresponding to $|\xi|^2 \pm i\Omega \eta$ are orthogonal to $\xi$. 
Since $\tilde u$ is divergence-free, we can forget about the
first eigenvector, and we obtain the representation formula
\begin{equation}\label{uformula}
  \tilde u_n(t,k) \,=\, e^{-t|\xi|^2}\Bigl(
  e^{-it\Omega\eta}\langle \tilde u_n^0(k),w_n^+(k)\rangle
  + e^{it\Omega\eta}\langle \tilde u_n^0(k),w_n^-(k)\rangle
\Bigr)~, \quad t \ge 0~,
\end{equation}
where $\tilde u_n^0(k) = \tilde u_n(0,k)$ and $\langle\cdot,\cdot
\rangle$ denotes the usual scalar product in $\CC^3$.

To estimate the norm of $\tilde u$ in the space 
$L^1(\RR_+,L^\infty(\DD))$, we proceed as in \cite{cdgg}. 
Using standard approximation arguments, it is easy to show that
\[
  \|\tilde u\|_{L^1(\RR_+,L^\infty(\DD))} \,=\,
  \sup_{\phi \in {\cal E}} <\tilde u,\phi>_{L^2(\RR_+,L^2(\DD))}~,
\]
where ${\cal E} = \{\phi \in C_0^\infty(\DD)\,|\,\|\phi\|_{L^{\infty}
(\RR_+,L^1(\DD))} \le 1 \}$. By the Parseval relation, we thus have
\begin{equation}\label{auxdis}
  \|\tilde u\|_{L^1(\RR_+,L^{\infty}(\DD))} \,=\, \sup_{\phi \in {\cal E}}
  \int_0^\infty\!\!\int_{\RR^2}\sum_{n \in \ZZ^*}\tilde u_n(t,k)
  \,\overline{\phi_n(t,k)} \dd k \dd t~,
\end{equation}
where $\phi_n(t,k)$ denotes of course the Fourier transform of
$\phi(t,\cdot)$. The idea is now to replace \eqref{uformula} into
\eqref{auxdis}, and to estimate the right-hand side. Before doing
that, we recall that the initial data $\tilde u_n^0(k)$ were 
assumed to vanish outside a finite ball ${\cal B}_R$ in Fourier 
space, see Proposition~\ref{strichartz}. In view of \eqref{uformula}, 
the same property holds for $\tilde u_n(t,k)$ for all $t \ge 0$. 
Thus $\tilde u_n(t,k)\equiv \psi_n(k)\tilde u_n(t,k)$, 
where
\begin{equation}\label{psindef}
  \psi_n(k) \,=\, (1-\delta_{n,0})\,\chi\Bigl(\frac{\sqrt{
  |k|^2+4\pi^2n^2}}{2R}\Bigr)~, \quad k \in \RR^2~, \quad 
  n \in \ZZ~.
\end{equation}
Here $\chi$ is as in \eqref{PRdef}, and $\delta$ is the Kronecker 
symbol. 

Given any $A \ge 0$ and any $B \in \RR$, we denote by $K[A,B] 
\in C^{\infty}(\DD)$ the function defined in Fourier variables 
by
\begin{equation}\label{Kdef}
  K[A,B]_n(k) \,=\, \frac{1}{2\pi}\,e^{-A |\xi|^2 + iB\eta}\,\psi_n(k)^2~, 
  \quad k \in \RR^2~, \quad n \in \ZZ~,
\end{equation}
where $|\xi|$ and $\eta$ are as in \eqref{xietadef}. The 
following estimate will be crucial: 

\begin{lm}\label{dispersion}
For any $R>0$ there exists $C_R > 0$ such that, for any 
$A \ge 0$ and any $B \in \RR$, the function $K[A,B] \in 
C^{\infty}(\DD)$ defined by \eqref{Kdef} satisfies
\[
  \|K[A,B]\|_{L^{\infty}(\DD)} \,\le\, C_R \frac{e^{-4\pi^2 A}}
  {\sqrt{|B|}}~.
\]
\end{lm}

We postpone the proof of this lemma and first conclude the 
proof of Proposition~\ref{strichartz}. After replacing \eqref{uformula} 
into \eqref{auxdis}, we have to estimate for each $\phi \in {\cal E}$ 
the quantity $M_+ + M_-$, where
\[
  M_\pm \,=\, \int_{\RR^2}\sum_{n \in \ZZ^*} \langle\tilde u_n^0(k),
  w_n^\pm(k)\rangle\, \Bigl\{\int_0^\infty e^{-t|\xi|^2\mp it\Omega\eta}
  \psi_n(k) \,\overline{\phi_n(t,k)}\dd t\Bigr\} \dd k~.
\]
Since the eigenvectors $w_n^\pm(k)$ are normalized, the Cauchy-Schwarz 
inequality and the Parseval relation imply that $|M_\pm| \le 
\|\tilde u_0\|_{L^2(\DD)} N_\pm$, where
\begin{align*}
  N_\pm^2 \,&=\, \int_{\RR^2}\sum_{n \in \ZZ^*}
  \Bigl|\int_0^\infty e^{-t|\xi|^2\mp it\Omega\eta}\,\psi_n(k)
  \,\overline{\phi_n(t,k)}\dd t\Bigr|^2\dd k \\ 
  \,&=\, \int_{\RR^2}\sum_{n \in \ZZ^*}\int_0^\infty\!\!\int_0^\infty 
  e^{-(t+s)|\xi|^2 \pm i(s-t)\Omega\eta}\,\psi_n(k)^2 
  \,\overline{\phi_n(t,k)}\,\phi_n(s,k)\dd t \dd s \dd k\\
  \,&=\, \int_0^\infty \!\!\int_0^\infty <K[t+s,\pm\Omega(s-t)]
  *\phi(s,\cdot)\,,\,\phi(t,\cdot)>_{L^2(\DD)}\dd t\dd s~.
\end{align*}
In the last line we have used the definition \eqref{Kdef} of 
$K[A,B]$ and the Parseval relation again. Now, since 
$\phi \in {\cal E}$, it follows from Young's inequality that
\begin{align*}
  |<K[A,B]*\phi(s,\cdot)\,,\,\phi(t,\cdot)>_{L^2(\DD)}|
  \,&\le\, \|K[A,B]\|_{L^\infty(\DD)}\,\|\phi(t,\cdot)\|_{L^1(\DD)} 
    \,\|\phi(s,\cdot)\|_{L^1(\DD)} \\
  \,&\le\, \|K[A,B]\|_{L^\infty(\DD)}~.
\end{align*}
Thus, setting $A = t+s$, $B = \pm\Omega (s-t)$, we obtain from 
Lemma~\ref{dispersion}
\[
  N_\pm^2 \,\le\, C_R \int_0^\infty \!\!\int_0^\infty 
  \frac{e^{-4\pi^2(t+s)}}{\sqrt{|\Omega||t-s|}}\dd t\dd s
  \,\le\, \frac{C_R}{|\Omega|^{1/2}}~.
\]
Summarizing, we have shown that $|M_\pm| \le C_R |\Omega|^{-1/4}
\|\tilde u_0\|_{L^2(\DD)}$ for all $\phi \in {\cal E}$, which in
turn implies $\|\tilde u\|_{L^1(\RR_+,L^\infty(\DD))} \le C_R 
|\Omega|^{-1/4}\|\tilde u_0\|_{L^2(\DD)}$. This concludes the 
proof of Proposition~\ref{strichartz}.  
\rule{2mm}{2mm}

\bigskip
\noindent{\bf Proof of lemma \ref{dispersion}:} Given $A \ge 0$ and
$B \in \RR$, we have to estimate the expression
\[
  K[A,B](x,z) \,=\, \frac{1}{4\pi^2}\int_{\RR^2} \sum_{n \in \ZZ^*} 
  e^{-A|\xi(k,n)|^2 +i B\eta(k,n)}\,\psi_n(k)^2 \,e^{i(k\cdot x 
  +2\pi nz)}\dd k~,
\]
where $|\xi|, \eta$ are defined in \eqref{xietadef} and $\psi_n(k)$
is given by \eqref{psindef}. Here again, we follow the approach 
presented in \cite[Chap.~5]{cdgg}. As $K[A,B](x,z)$ is a radially 
symmetric function of $x \in \RR^2$, we can assume without loss 
of generality that $x_2 = 0$. Clearly, we can also suppose that 
$B \ge 0$. Let $L$ be the first-order differential operator defined by
\[
  L \,=\, \frac{1}{1+B\alpha(k,n)^2}(1 + i\alpha(k,n)
  \partial_{k_2})~, \quad \hbox{where}\quad \alpha(k,n) = 
  -\partial_{k_2}\eta(k,n)~.
\]
Then $L(e^{iB\eta(k,n)}) = e^{iB\eta(k,n)}$, and integrating by 
parts (over the variable $k \in \RR^2$) we find
\[
  K[A,B]((x_1,0),z) \,=\, \frac{1}{4\pi^2}\int_{\RR^2} \sum_{n \in \ZZ^*} 
  e^{iB\eta(k,n)} \,e^{i(k_1x_1 +2\pi nz)} L^t\Bigl(e^{-A|\xi(k,n)|^2} 
  \,\psi_n(k)^2\Bigr)\dd k~,
\]
where $L^t$ denotes the formal adjoint of $L$. A direct calculation
gives
\begin{align*}
  L^t\Bigl(e^{-A|\xi(k,n)|^2} \,\psi_n(k)^2\Bigr) \,=\, 
  &\Bigl(\frac{1}{1+B\alpha^2} -i(\partial_{k_2}\alpha)
  \frac{1-B\alpha^2}{(1+B\alpha^2)^2}\Bigr)e^{-A|\xi(k,n)|^2}
  \,\psi_n(k)^2 \\
  \,&-\, \frac{i\alpha}{1+B\alpha^2}\,\partial_{k_2}
  \Bigl(e^{-A|\xi(k,n)|^2} \,\psi_n(k)^2\Bigr)~.
\end{align*}
We have to estimate this quantity for $(k,n) \in {\cal B}_{2R}$ and $n \neq 0$, 
because $\psi_n(k) = 0$ if $(k,n) \notin {\cal B}_{2R}$ or $n = 0$. We first 
observe that
\[
  |\xi(k,n)| \,\ge\, 2\pi~, \quad \hbox{and}\quad
  |\alpha(k,n)| \,=\, \frac{2\pi|n| |k_2|}{(|k|^2+4\pi^2n^2)^{3/2}}
  \,\ge\, \frac{\pi|k_2|}{4R^3}~.
\]
Moreover, there exists $C_R > 0$ such that $|\alpha(k,n)| + 
|\partial_{k_2}\alpha(k,n)| \le C_R$. As a consequence, we have
\[
  \frac{1}{1+B\alpha^2} + \frac{|1-B\alpha^2|}{(1+B\alpha^2)^2}
  +\frac{|\alpha|}{1+B\alpha^2} \,\le\, \frac{C_R}{1+Bk_2^2}~,
\]
so that $|L^t(e^{-A|\xi(k,n)|^2} \,\psi_n(k)^2)| \le C_R\,e^{-4\pi^2A}
\psi_n(k)(1+Bk_2^2)^{-1}$. We conclude that
\begin{align*}
  \|K[A,B]\|_{L^\infty(\DD)} \,&\le\, \frac{1}{4\pi^2}\int_{\RR^2} 
  \sum_{n \in \ZZ^*} |L^t(e^{-A|\xi(k,n)|^2} \,\psi_n(k)^2)| \dd k \\
   \,&\le\, C_R\,e^{-4\pi^2A} \int_\RR \frac{\dd k_2}{1+Bk_2^2}
  \,\le\, C_R \,\frac{e^{-4\pi^2A}}{\sqrt{B}}~,
\end{align*}
which is the desired estimate. \rule{2mm}{2mm}


\end{document}